\documentclass[a4paper,11pt]{article}

\usepackage[leqno]{amsmath}
\usepackage{amssymb,amsthm,ifthen}
\usepackage{url,enumitem}
\usepackage[pdftex]{graphicx}
\usepackage[curve,knot]{xypic}
\usepackage{color}

\usepackage{tikz}
\usetikzlibrary{calc}
\usepgflibrary{shapes.geometric}
\usepgflibrary{shapes.misc}
\usetikzlibrary{positioning}
\usetikzlibrary{decorations}

\usepackage[a4paper,left=1in,width=458pt,textheight=650pt,top=1.5in]{geometry}

\usepackage{bbm}

\DeclareMathAlphabet{\mathitbf}{OML}{cmm}{b}{it}

\usepackage{hyperref}

\input cyracc.def

\providecommand{\keywords}[1]{\noindent \footnotesize {\textbf{Keywords:} #1 } \normalsize }

\newboolean{Sectext}
\newsavebox{\MSCtextP}
\newsavebox{\MSCtextS}
\newcommand{\MSC}[2]{
\sbox{\MSCtextP}{#1}
\ifthenelse{\equal{#2}{NONE}}{\setboolean{Sectext}{false}}{\setboolean{Sectext}{true}}
\sbox{\MSCtextS}{#2}
}
\newcommand{\printMSC}{
\ifthenelse{\boolean{Sectext}}{\usebox{\MSCtextP} (Primary), \usebox{\MSCtextS}}{\usebox{\MSCtextP}}
}
\newcommand{\printMSCwithtitle}{
\ifthenelse{\boolean{Sectext}}{\noindent \footnotesize \textbf{Mathematics Subject Classification (2000):} \usebox{\MSCtextP} (Primary), \usebox{\MSCtextS} \normalsize}{\noindent \footnotesize \textbf{Mathematics Subject Classification (2000):} \usebox{\MSCtextP} \normalsize}
}

\hypersetup{
pdftitle={Quantum cluster algebra structures on quantum Grassmannians and their quantum Schubert cells: the finite-type cases},
pdfauthor={Jan E. Grabowski and Stephane Launois},
pdfstartview=FitH
}

\newcommand{\balpha}{\boldsymbol\alpha}

\newcommand{\bbeta}{\boldsymbol\beta}

\newcommand{\card}[1]{\lvert #1 \rvert}

\let\chisave\chi
\renewcommand{\chi}{{%
 \mathchoice{\raisebox{0.25ex}{$\displaystyle\chisave$}}
            {\raisebox{0.2ex}{$\textstyle\chisave$}}
            {\raisebox{0.2ex}{$\scriptstyle\chisave$}}
            {\raisebox{0.1ex}{$\scriptscriptstyle\chisave$}}}}

\newcommand{\complex}{\ensuremath \mathbb{C}}

\newcommand{\cross}{\times}
\newcommand{\curly}[1]{\ensuremath{\mathcal{#1}}}

\newcommand{\defeq}{\stackrel{\scriptscriptstyle{\mathrm{def}}}{=}}

\newcommand{\dual}[1]{\ensuremath {#1}^{*}}

\renewcommand{\epsilon}{\varepsilon}

\newcommand{\inj}{\hookrightarrow}

\newcommand{\iso}{\ensuremath \cong}
\newcommand{\itbf}[1]{\mathitbf{#1}}

\newcommand{\matentrybf}[1]{\mathitbf{X}_{\mathbf{#1}}}

\newcommand{\nat}{\ensuremath \mathbb{N}}

\renewcommand{\phi}{\varphi}

\newcommand{\qminor}[1]{\Delta_{q}^{#1}}
\newcommand{\qminorbf}[1]{\mathbf{\Delta}_{\mathitbf{q}}^{\mathbf{#1}}}

\newcommand{\ie}{i.e.\ }
\newcommand{\eg}{e.g.\ }

\theoremstyle{plain}
\newtheorem{theorem}{Theorem}
\newtheorem*{theorem*}{Theorem}
\newtheorem{proposition}[theorem]{Proposition}

\theoremstyle{definition}
\newtheorem{definition+}[theorem]{Definition}
\theoremstyle{remark}
\newtheorem{remark}[theorem]{Remark}
\newtheorem*{example*}{Example}
\newtheorem*{examplectd*}{Example (continued)}

\title{Quantum cluster algebra structures on quantum Grassmannians \\ and their quantum Schubert cells: the finite-type cases}
\author{Jan E. Grabowski\footnotemark[2] 
\\ \small{\textit{Mathematical Institute, University of Oxford,}}
\\ \small{\textit{24-29 St.\ Giles', Oxford, OX1 3LB, United Kingdom}}
\and St\'{e}phane Launois\footnotemark[3]
\\ \small{\textit{School of Mathematics, Statistics and Actuarial Science,
University of Kent,}}
\\ \small{\textit{Canterbury, CT2 7NF, United Kingdom}}
}
\date{22nd December 2009}

\begin{document}

\maketitle

\renewcommand{\thefootnote}{\fnsymbol{footnote}}
\footnotetext[2]{Email: \url{jan.grabowski@maths.ox.ac.uk}.  Website: \url{http://people.maths.ox.ac.uk/~grabowsk/}}
\footnotetext[3]{Email: \url{S.Launois@kent.ac.uk}.  Website: \url{http://www.kent.ac.uk/ims/personal/sl261/}}
\renewcommand{\thefootnote}{\arabic{footnote}}
\setcounter{footnote}{0}

\begin{abstract} 
\noindent We exhibit quantum cluster algebra structures on quantum Grassmannians $\mathbb{K}_{q}[\mathrm{Gr}(2,n)]$ and their quantum Schubert cells, as well as on $\mathbb{K}_{q}[\mathrm{Gr}(3,6)]$, $\mathbb{K}_{q}[\mathrm{Gr}(3,7)]$ and $\mathbb{K}_{q}[\mathrm{Gr}(3,8)]$.  These cases are precisely those where the quantum cluster algebra is of finite type and the structures we describe quantize those found by Scott for the classical situation.
\end{abstract}

\MSC{\footnotesize 20G42}{\footnotesize 16W35, 17B37}
\printMSCwithtitle

\vfill 

\tableofcontents
\vfill

\pagebreak

\section{Introduction}\label{s:intro}

Because of their wide range of connections with other areas of Mathematics, cluster algebras have been extensively studied in the recent years and many well-known algebras have been shown to have a cluster algebra structure. For instance, Scott (\cite{Scott-Grassmannians}) proved that the homogeneous coordinate ring of the Grassmannian $\mathrm{Gr}(k,n)$ provides an example of a cluster algebra. Among these examples, only a few are of finite type (i.e. have only finitely many cluster variables): 
$\mathbb{C}[\mathrm{Gr}(2,n)] $ which is of type $A_{n-3}$, 
$\mathbb{C}[\mathrm{Gr}(3,6)] $ which is of type $D_{4}$, 
$\mathbb{C}[\mathrm{Gr}(3,7)] $ which is of type $E_{6}$ and 
$\mathbb{C}[\mathrm{Gr}(3,8)] $ which is of type $E_{8}$. 
Scott's result was then generalised by Gei\ss, Leclerc and Schr\"{o}er (\cite{GLS-PFV}) who have shown that the multi-homogeneous coordinate rings of partial flag varieties and their associated unipotent radicals have a cluster algebra structure.  

On the contrary, very few examples of quantum cluster algebras are known.  Quantum cluster algebras were introduced and studied by Berenstein and Zelevinsky in \cite{BZ-QCA}. In particular, in this paper they conjecture that the quantized coordinate ring of a double Bruhat cell is a quantum cluster algebra. To the best of our knowledge this conjecture remains unproven, whereas the classical counterpart of this conjecture was proved by Berenstein, Fomin and Zelevinsky in \cite{BFZ-CA3}. 

In \cite{QCAexamples}, the first named author proved that the quantum Grassmannian $\mathbb{K}_q[\mathrm{Gr}(2,5)] $ is a quantum cluster algebra of type $A_2$. In the present paper, we extend this example and show that the quantized coordinate rings $\mathbb{K}_{q}[\mathrm{Gr}(2,n)]$ for $n\geq 3$ are quantum cluster algebras of type $A_{n-3}$ (Section~\ref{ss:Gr2n}).  By computer-aided calculation, we also show that $\mathbb{K}_{q}[\mathrm{Gr}(3,6)]$, $\mathbb{K}_{q}[\mathrm{Gr}(3,7)]$ and $\mathbb{K}_{q}[\mathrm{Gr}(3,8)]$ are quantum cluster algebras of type $D_4$, $E_6$ and $E_8$ respectively (Section~\ref{ss:Gr3-678}).

From this, we have also obtained quantum cluster algebra structures on the quantum Schubert cells of the $k=2$ Grassmannians (Section~\ref{s:qSchubertCells}).  The quantum Schubert cell associated to the partition $(t,s)$ (where $t\geq s$ and  $t,s\leq n-2$) is of quantum cluster algebra type $A_{s-1}$, independent of $t$. 

We view these results as a step towards achieving the goal of quantizing the work of Gei\ss, Leclerc and Schr\"{o}er (\cite{GLS-PFV}). 

In order to prove these results, we proceed as follows. First, in each case, we define an appropriate quantum initial seed. Roughly speaking, a quantum initial seed is formed of two components: a quiver and a matrix. The quiver is easy to define in our cases: we define it in identical fashion to the classical case.  In particular, it will always be mutation-equivalent to a Dynkin quiver (since we are only considering the Grassmannians $\mathrm{Gr}(2,n)$ for $n \geq 4$ and $\mathrm{Gr}(3,n)$ for $n \in \{6,7,8\}$), so that the quantum cluster algebra associated is of finite type. The matrix is also easy to define: it encodes the quasi-commutation relations between the elements of the quantum initial seed. One of the main difficulties here is to check that this matrix satisfies a certain compatibility condition with the quiver. Once we know that this compatibility condition is satisfied, we then identify explicitly all the quantum cluster variables (there are only finitely many of them), and prove our main results.  

\subsection*{Acknowledgements}

The first named author would like to acknowledge the provision of facilities by Keble College and the Mathematical Institute in Oxford. The research of the second named author was supported by a Marie Curie European Reintegration Grant within the $7^{\text{th}}$ European Community Framework Programme.

\section{Recollections on cluster algebras and their quantum analogues}\label{s:CAsandQCAs}

\subsection{Cluster algebras}\label{ss:CAs}

We will recall briefly the definition of a cluster algebra of geometric type (\cite{FZ-CA1}) and describe in detail only the ``no coefficients'' case.  We start with an \emph{initial seed} $(\underline{y},B)$, consisting of a tuple of generators (called a \emph{cluster}) for the cluster algebra and an \emph{exchange matrix} $B=(b_{ij})$ with integral entries.  (A cluster is not a complete set of generators, but a subset of such a set.)  More seeds are obtained via \emph{mutation} of the initial seed.  Matrix mutation $\mu_{k}$ is involutive and given by the rule
\[ (\mu_{k}(B))_{ij} = \begin{cases} -b_{ij} & \text{if}\ i=k\ \text{or}\ j=k \\ b_{ij}+\frac{|b_{ik}|b_{kj}+b_{ik}|b_{kj}|}{2} & \text{otherwise.} \end{cases} \]
For example,
\begin{align*}
\begin{pmatrix} 0 & 1 & 0 \\ -1 & 0 & 1 \\ 0 & -1 & 0 \end{pmatrix} & \stackrel{\mu_{1}}{\longrightarrow} \begin{pmatrix} 0 & -1 & 0 \\ 1 & 0 & 1 \\ 0 & -1 & 0 \end{pmatrix} & \begin{pmatrix} 0 & 1 & 0 \\ -1 & 0 & 1 \\ 0 & -1 & 0 \end{pmatrix} & \stackrel{\mu_{2}}{\longrightarrow} \begin{pmatrix} 0 & -1 & 1 \\ 1 & 0 & -1 \\ -1 & 1 & 0 \end{pmatrix}.
\end{align*}

If $(\underline{y}=(y_{1},\ldots ,y_{d}),B)$ is the initial seed then the mutated seed in direction $k$ is given by $(\mu_{k}(\underline{y})=(y_{1},\ldots,y_{k}^{\ast},\ldots,y_{d}),\mu_{k}(B))$, where the new generator $y_{k}^{\ast}$ is determined by the \emph{exchange relation}
\[ y_{k}y_{k}^{\ast}=\prod_{b_{ik}>0} y_{i}^{b_{ik}} + \prod_{b_{ik}<0} y_{i}^{-b_{ik}}. \]

The alternative quiver description converts $B$ to a quiver by the rule that a strictly positive entry $b_{ij}$ determines a weighted arrow $i\stackrel{b_{ij}}{\to} j$ and a strictly negative one a weighted arrow in the opposite direction.  (Thus $B$ is what is sometimes termed an incidence matrix for the quiver; the adjacency matrix is the matrix obtained from the incidence matrix by replacing $b_{ij}$ by $\max\{0,b_{ij}\}$.)  Then matrix mutation defines the operation of quiver mutation, for example
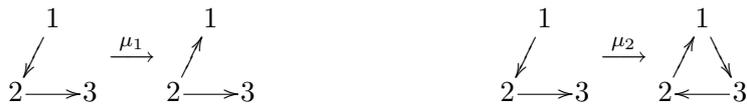
\begin{align*}
\raisebox{1.5em}[3em][2em]{\xymatrix@1@=1pt@!{ {} & {1} \ar[ddl] & {} \\ {} & {} & {} \\ {2} \ar[rr] & {} & {3} }} & \stackrel{\mu_{1}}{\longrightarrow} \raisebox{1.5em}[0em][2em]{\xymatrix@1@=1pt@!{ {} & {1} & {} \\ {} & {} & {} \\ {2} \ar[rr] \ar[uur] & {} & {3} }} & \raisebox{1.5em}[0em][2em]{\xymatrix@1@=1pt@!{ {} & {1} \ar[ddl] & {} \\ {} & {} & {} \\ {2} \ar[rr] & {} & {3} }} & \stackrel{\mu_{2}}{\longrightarrow} \raisebox{1.5em}[0em][2em]{\xymatrix@1@=1pt@!{ {} & {1} \ar[ddr] & {} \\ {} & {} & {} \\ {2} \ar[uur] & {} & {3} \ar[ll] }} 
\end{align*}

We say an algebra $\curly{A}$ is a cluster algebra or admits a cluster algebra structure if the set of all cluster variables (\ie the union of all the clusters) is a generating set for $\curly{A}$.  A cluster algebra is of finite type (as all our examples will be) if the quiver of $B$ lies in the same mutation equivalence class as an orientation of a finite-type Dynkin diagram and the type of the cluster algebra is the type of this diagram.

The ``with coefficients'' version includes additional generators present in every cluster that are never mutated but monomials in them also appear as coefficients in the exchange relations.  In the quiver approach, these correspond to ``frozen'' vertices, indicated by drawing a box around the vertex.  We will refer to the elements of clusters that are not coefficients as \emph{mutable} cluster variables.  We will indicate the mutable variables in a cluster by boldface type.  If the cluster algebra under consideration is of finite type $X_{n}$, there is a bijection between the set of all mutable cluster variables (from all clusters) and the almost positive roots of the root system of type $X_{n}$.  (The almost positive roots are the positive roots together with the negative simple roots.)

\subsection{Quantum cluster algebras}\label{ss:QCAs}

Berenstein and Zelevinsky (\cite{BZ-QCA}) have given a definition of a quantum cluster algebra.  These algebras are now non-commutative but not so far from being commutative.  Each quantum seed includes an additional piece of data, an integral skew-symmetric matrix $L=(l_{ij})$ describing \emph{quasi-commutation} relations between the variables in the cluster.  Quasi-commuting means $ab=q^{l_{ab}}ba$, also written in $q$-commutator notation as $[a,b]_{q^{l_{ab}}}=0$.  

There is also a mutation rule for these quasi-commutation matrices and a modified exchange relation that involves further coefficients that are powers of $q$ derived from $B$ and $L$, which we describe now.  Given a quantum cluster $\underline{y}=(X_{1},\ldots,X_{r})$, exchange matrix $B$ and quasi-commutation matrix $L$, the exchange relation for mutation in the direction $k$ is given by
\[ X_{k}'=M(-\underline{\boldsymbol{e}}_{k}+\sum_{b_{ik}>0}b_{ik}\underline{\boldsymbol{e}}_{i})+M(-\underline{\boldsymbol{e}}_{k}-\sum_{b_{ik}<0}b_{ik}\underline{\boldsymbol{e}}_{i}) \]
where the vector $\underline{\boldsymbol{e}}_{i}\in \complex^{r}$ ($r$ being the number of rows of $B$) is the $i$th standard basis vector and 
\[ M(a_{1},\dotsc ,a_{r})=q^{\frac{1}{2}\sum_{i<j} a_{i}a_{j}l_{ji}}X_{1}^{a_{1}}\dotsm X_{r}^{a_{r}}. \]
By construction, the integers $a_{i}$ are all non-negative except for $a_{k}=-1$.  The monomial $M$ (as we have defined it here) is related to the concept of a toric frame, also introduced in \cite{BZ-QCA}.  The latter is a technical device used to make the general definition of a quantum cluster algebra.  For our examples, where we start with a known algebra and want to exhibit a quantum cluster algebra structure on this, it will suffice to think of $M$ simply as a rule determining the exchange monomials.  

We note that the presence of the ``$\frac{1}{2}$'' factor in the definition of $M$ suggests that if we wanted to work over fields other than $\complex$, we may need to extend scalars by introducing a square root of $q$.  In fact this will not be necessary in all examples but it would be required in some. 

The natural requirement that all mutated clusters also quasi-commute leads to a compatibility condition between $B$ and $L$, namely that $B^{T}L$ consists of two blocks, one diagonal with positive integer diagonal entries and one zero.  (However, these blocks need not be contiguous, depending on the ordering of the row and column labels.)  We will denote by $0_{m,n}$ the $m\cross n$ zero matrix and by $I_{m}$ the $m\cross m$ identity matrix.  Block matrices will be written in the usual way, \eg $(A\ B\ C)$ or $\left( \begin{smallmatrix} A & B \\ C & D \end{smallmatrix} \right)$.

Importantly, Berenstein and Zelevinsky show that the exchange graph (whose vertices are the clusters and edges are mutations) remains unchanged in the quantum setting.  That is, the matrix $L$ does not influence the exchange graph.  It follows that quantum cluster algebras are classified by Dynkin types in exactly the same way as the classical cluster algebras.

Other previously-known examples of quantum cluster algebras include quantum symmetric algebras (of rank 0) and conjecturally quantum double Bruhat cells (\cite{BZ-QCA}). 

\section{Quantum Grassmannians}\label{s:qGrassmannians}

Throughout, $\mathbb{C}$ denotes the field of complex numbers and $\mathbb{K}$ is a field. Moreover we assume that $q \in \mathbb{K} $ is a non-zero element such that $q^{1/2}$ exists in $\mathbb{K}$. Let $C$ be an $l\cross l$ generalized Cartan matrix with columns indexed by a set $I$.  Let $(H,\Pi,\Pi^{\vee})$ be a minimal realization of $C$, where $H \iso \complex^{2\left| I \right|-\text{rank}(C)}$, $\Pi=\{ \alpha_{i} \mid i \in I \} \subset \dual{H}$ (the simple roots) and $\Pi^{\vee}=\{ h_{i} \mid i \in I \} \subset H$ (the simple coroots).  Then we say $\curly{C}=(C,I,H,\Pi,\Pi^{\vee})$ is a root datum associated to $C$.  (Lusztig (\cite{LusztigBook}) has a more general definition of a root datum but this one will suffice for our purposes.)

If $G=G(\curly{C})$ is a connected semisimple complex algebraic group associated to $\curly{C}$, $G$ has a (standard) parabolic subgroup $P_{J}$ associated to any choice of subset $J\subseteq I$.  From this we can form $G/P_{J}$, a partial flag variety; the choice $J=\emptyset$ gives $G/P_{\emptyset}=G/B$, the full flag variety.  We set $D=I\setminus J$.

The partial flag variety $G/P_{J}$ is a projective variety, via the well-known Pl\"{u}cker embedding $G/P_{J} \inj \prod_{d\in D} \mathbb{P}(L(\omega_{d}))$.  (Here, $L(\lambda)$ is the irreducible $G$-module corresponding to a dominant integral weight $\lambda$ and $\{\omega_{i}\}_{i\in I}$ are the fundamental weights.)  Via the Pl\"{u}cker embedding, we may form the corresponding $\nat^{D}$-graded multi-homogeneous coordinate algebra $\complex[G/P_{J}]=\bigoplus_{\lambda \in \nat^{D}} \dual{L(\lambda)}$.

The coordinate ring $\complex[G]$ has a quantum analogue, $\mathbb{K}_{q}[G]$ (see for example \cite{Brown-Goodearl}, where this algebra is denoted $\curly{O}_{q}(G)$).  Via this quantized coordinate ring, we can define a quantization of $\mathbb{K}[G/P_{J}]$, $\mathbb{K}_{q}[G/P_{J}]$.

The case we consider is that of the partial flag variety obtained from $G=G(A_{n})=SL_{n+1}(\complex)$ and $J=I \setminus \{k\}$, namely $G/P_{J}=\mathrm{Gr}(k,n)$, the Grassmannian of $k$-dimensional subspaces in $\complex^{n}$. In this case, the quantized coordinate ring $\mathbb{K}_{q}[\mathrm{Gr}(k,n)]$ is the subalgebra of the quantum matrix algebra $\mathbb{K}_{q}[M(k,n)]$ generated by the quantum Pl\"{u}cker coordinates.

We recall that the quantum matrix algebra $\mathbb{K}_{q}[\mathrm{M}(k,n)]$ is the $\mathbb{K}$-algebra generated by the set $\{ X_{ij} \mid 1\leq i\leq k,\ 1\leq j \leq n \}$ subject to the quantum $2\cross 2$ matrix relations on each $2\cross 2$ submatrix of \[ \begin{pmatrix} X_{11} & X_{12} & \cdots & X_{1n} \\ \vdots & \vdots & \ddots & \vdots \\ X_{k1} & X_{k2} & \cdots & X_{kn} \end{pmatrix}, \] where the quantum $2\cross 2$ matrix relations on $\left( \begin{smallmatrix} a & b \\ c & d \end{smallmatrix} \right)$ are
\begin{align*} ab & = qba & ac & = qca & bc & = cb \\ bd & = qdb & cd & = qdc & ad-da & = (q-q^{-1})bc. \end{align*}

 The quantized coordinate ring $\mathbb{K}_{q}[\mathrm{Gr}(k,n)]$ is the subalgebra of the quantum matrix algebra $\mathbb{K}_{q}[\mathrm{M}(k,n)]$ generated by the quantum Pl\"{u}cker coordinates. In other words,  $\mathbb{K}_{q}[\mathrm{Gr}(k,n)]$ is the subalgebra of $\mathbb{K}_{q}[\mathrm{M}(k,n)]$ generated by the $k\cross k$ quantum minors of $\mathbb{K}_{q}[\mathrm{M}(k,n)]$. 

Recall that the $k\cross k$ quantum minor $\Delta_{q}^{I}$ associated to the $k$-subset $I=\{ i_{1} < i_{2} < \cdots < i_{k} \}$ of $\{1, \dots , n\}$ is defined to be
\[ \Delta_{q}^{I} \defeq \sum_{\sigma \in S_{k}} (-q)^{l(\sigma)}X_{1i_{\sigma(1)}}\cdots X_{k\mspace{0.5mu}i_{\sigma(k)}} \] where $S_{k}$ is the symmetric group of degree $k$ and $l$ is the usual length function on this.  (In fact we are considering the quantum minor $\Delta_{\{1,\ldots,k\}}^{I}$ but since we are working in $\mathbb{K}_{q}[\mathrm{M}(k,n)]$ there is no other choice for the row-minor subset and so we omit it.  Quantum minors of smaller degree can be defined in the obvious way but we will not need these.)  For example, when $k=2$, the quantum minor $\Delta_{q}^{ij}$ for $i<j$ is equal to $X_{1i}X_{2j}-qX_{1j}X_{2i}$.

Then we denote by $\curly{P}_{q}$ the set of all quantum Pl\"{u}cker coordinates, that is $$\curly{P}_{q} = \{ \Delta_{q}^{I}  \mid I \subseteq \{1,\ldots ,n\}, \card{I}=k \}.$$
This is a generating set of $\mathbb{K}_{q}[\mathrm{Gr}(k,n)]$.
  
\subsection{The quantum Grassmannians \texorpdfstring{$\mathbb{K}_{q}[\mathrm{Gr}(2,n)]$}{Kq[Gr(2,n)]}}\label{ss:Gr2n}

The first case we consider is that of the partial flag variety obtained from $G=G(A_{n})=SL_{n+1}(\complex)$ and $J=I \setminus \{2\}$, namely $G/P_{J}=\mathrm{Gr}(2,n)$, the Grassmannian of 2-dimensional subspaces in $\complex^{n}$.  We give an initial quantum seed for a quantum cluster algebra structure on $\mathbb{K}_{q}[\mathrm{Gr}(2,n)]$.  For the initial quantum cluster we choose
\[ \underline{\tilde{y}} = ( \qminor{1n},\qminorbf{1(\itbf{n}-1)},\qminorbf{1(\itbf{n}-2)},\ldots,\qminorbf{14},\qminorbf{13},\qminor{12},\qminor{23},\qminor{34},\ldots,\qminor{(n-2)(n-1)},\qminor{(n-1)n}). \]
As described in \cite{Leclerc-Zelevinsky} and \cite{Scott-QMinors}, two quantum Pl\"{u}cker coordinates $\Delta_{q}^{ij}$ and $\Delta_{q}^{kl}$ quasi-commute when $\{i,j\}$ and $\{k,l\}$ are weakly separated, meaning---in this particular case---that the corresponding diagonals of a regular $n$-gon do not cross.  (The power of $q$ appearing in the corresponding quasi-commutation relation is also combinatorially determined.)  So, the above cluster is a set of quasi-commuting variables: the corresponding diagonals of the $n$-gon are seen to be the $n$ edges (in bijection with the coefficients) and $n-3$ non-crossing diagonals, $(1,i)$ for $3\leq i\leq n-1$.  That is, this cluster corresponds to a triangulation of the $n$-gon, as in the classical case (see \eg \cite{FZ-CA2}).

It is well known that $\mathbb{K}_{q}[\mathrm{Gr}(2,n)]$ is a noetherian domain, and so admits a skew-field of fractions that we denote by $\mathcal{F}_q$. Moreover the elements of the initial quantum cluster $\underline{\tilde{y}}$ generate $\mathcal{F}_q$ (as a skew-field). In what follows, all computations take place in $\mathcal{F}_q$. We now return to the quantum initial seed.

The corresponding quantum exchange matrix $B$ is equal to that for the well-known cluster algebra structure on $\complex[\mathrm{Gr}(2,n)]$ (\cite{Scott-Grassmannians}) and, along with its quiver $\Gamma(B)$, is as follows.  The matrix $B$ has one row for each entry of $\underline{\tilde{y}}$ (in that order) and has columns indexed by the mutable cluster variables, \ie $(\qminorbf{1(\itbf{n}-1)},\qminorbf{1(\itbf{n}-2)},\ldots,\qminorbf{13})$.  For brevity, we use just the minor label, which we will write $[ij]$, rather than $\qminor{ij}$.  We describe the column of $B$ indexed by $\qminorbf{1\itbf{k}}$ ($3\leq k \leq n-1$):
\begin{align*}
B_{[1i][1k]} & =  
  \begin{cases} 
    0 & \text{for}\ i\geq k+2 \\
   -1 & \text{for}\ i=k+1 \\
    0 & \text{for}\ i=k \\
    1 & \text{for}\ i=k-1 \\
    0 & \text{for}\ i\leq k-2
  \end{cases} &
B_{[j(j+1)][1k]} & =
  \begin{cases} 
    0 & \text{for}\ j\leq k-2 \\
   -1 & \text{for}\ j=k-1 \\
    1 & \text{for}\ j=k \\
    0 & \text{for}\ j\geq k+1
  \end{cases}
\end{align*}
where $2\leq i\leq n$ and $2\leq j \leq n-1$.  For example, for $n=8$ we have 
\[ B = 
\bordermatrix{
 & \scriptstyle{[17]} & \scriptstyle{[16]} & \scriptstyle{[15]} & \scriptstyle{[14]} & \scriptstyle{[13]} \cr
\scriptstyle{[18]} & -1 & 0 & 0 & 0 & 0 \cr
\scriptstyle{[17]} & \mathbf{0} & \mathbf{-1} & \mathbf{0} & \mathbf{0} & \mathbf{0} \cr
\scriptstyle{[16]} & \mathbf{1} & \mathbf{0} & \mathbf{-1} & \mathbf{0} & \mathbf{0}  \cr
\scriptstyle{[15]} & \mathbf{0} & \mathbf{1} & \mathbf{0} & \mathbf{-1} & \mathbf{0}  \cr
\scriptstyle{[14]} & \mathbf{0} & \mathbf{0} & \mathbf{1} & \mathbf{0} & \mathbf{-1}  \cr
\scriptstyle{[13]} & \mathbf{0} & \mathbf{0} & \mathbf{0} & \mathbf{1} & \mathbf{0}  \cr
\scriptstyle{[12]} & 0 & 0 & 0 & 0 & 1 \cr
\scriptstyle{[23]} & 0 & 0 & 0 & 0 & -1 \cr
\scriptstyle{[34]} & 0 & 0 & 0 & -1 & 1 \cr
\scriptstyle{[45]} & 0 & 0 & -1 & 1 & 0 \cr
\scriptstyle{[56]} & 0 & -1 & 1 & 0 & 0 \cr
\scriptstyle{[67]} & -1 & 1 & 0 & 0 & 0 \cr
\scriptstyle{[78]} & 1 & 0 & 0 & 0 & 0 
} \]
We note that $B$ has a natural block structure, as the submatrix on the row set $\{[1n],\ldots,[12]\}$ and that on the row set $\{[23],\ldots,[(n-1)n]\}$.

The corresponding quiver (for $n=8$), firstly with minor labels and secondly with diagonals of an octagon, is 
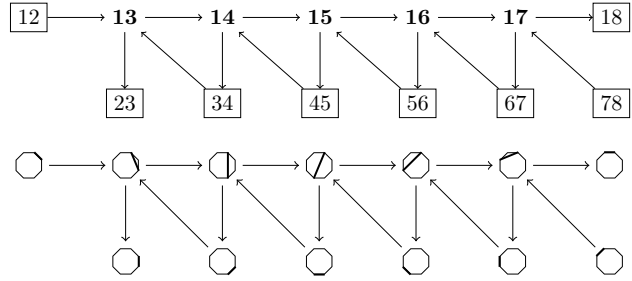
\begin{center} 
\scalebox{0.75}{\begin{tikzpicture}

\node (12) at (0,2) [rectangle,draw=black] {12};
\node (13) [right=of 12] {\textbf{13}};
\node (14) [right=of 13] {\textbf{14}};
\node (15) [right=of 14] {\textbf{15}};
\node (16) [right=of 15] {\textbf{16}};
\node (17) [right=of 16] {\textbf{17}};
\node (18) [rectangle,draw=black] [right=of 17] {18};
\node (23) [rectangle,draw=black] [below=of 13] {23};
\node (34) [rectangle,draw=black] [below=of 14] {34};
\node (45) [rectangle,draw=black] [below=of 15] {45};
\node (56) [rectangle,draw=black] [below=of 16] {56};
\node (67) [rectangle,draw=black] [below=of 17] {67};
\node (78) [rectangle,draw=black] [below=of 18] {78};

\draw[->] (12) to (13);
\draw[->] (13) to (14);
\draw[->] (13) to (23);
\draw[->] (14) to (15);
\draw[->] (14) to (34);
\draw[->] (15) to (16);
\draw[->] (15) to (45);
\draw[->] (16) to (17);
\draw[->] (16) to (56);
\draw[->] (17) to (18);
\draw[->] (17) to (67);

\draw[->] (34) to (13);
\draw[->] (45) to (14);
\draw[->] (56) to (15);
\draw[->] (67) to (16);
\draw[->] (78) to (17);

\end{tikzpicture}}
\end{center}

\begin{center}
\scalebox{0.925}{\begin{tikzpicture}

\node (Octagon12) at (0,2) [regular polygon, regular polygon sides=8, draw] {};
\draw[thick] (Octagon12.corner 1) -- (Octagon12.corner 8);

\node (Octagon13) [right=of Octagon12] [regular polygon, regular polygon sides=8, draw] {};
\draw[thick] (Octagon13.corner 1) -- (Octagon13.corner 7);

\node (Octagon14) [right=of Octagon13] [regular polygon, regular polygon sides=8, draw] {};
\draw[thick] (Octagon14.corner 1) -- (Octagon14.corner 6);

\node (Octagon15) [right=of Octagon14] [regular polygon, regular polygon sides=8, draw] {};
\draw[thick] (Octagon15.corner 1) -- (Octagon15.corner 5);

\node (Octagon16) [right=of Octagon15] [regular polygon, regular polygon sides=8, draw] {};
\draw[thick] (Octagon16.corner 1) -- (Octagon16.corner 4);

\node (Octagon17) [right=of Octagon16] [regular polygon, regular polygon sides=8, draw] {};
\draw[thick] (Octagon17.corner 1) -- (Octagon17.corner 3);

\node (Octagon18) [right=of Octagon17] [regular polygon, regular polygon sides=8, draw] {};
\draw[thick] (Octagon18.corner 1) -- (Octagon18.corner 2);

\node (Octagon23) [below=of Octagon13] [regular polygon, regular polygon sides=8, draw] {};
\draw[thick] (Octagon23.corner 8) -- (Octagon23.corner 7);

\node (Octagon34) [below=of Octagon14] [regular polygon, regular polygon sides=8, draw] {};
\draw[thick] (Octagon34.corner 7) -- (Octagon34.corner 6);

\node (Octagon45) [below=of Octagon15] [regular polygon, regular polygon sides=8, draw] {};
\draw[thick] (Octagon45.corner 6) -- (Octagon45.corner 5);

\node (Octagon56) [below=of Octagon16] [regular polygon, regular polygon sides=8, draw] {};
\draw[thick] (Octagon56.corner 5) -- (Octagon56.corner 4);

\node (Octagon67) [below=of Octagon17] [regular polygon, regular polygon sides=8, draw] {};
\draw[thick] (Octagon67.corner 4) -- (Octagon67.corner 3);

\node (Octagon78) [below=of Octagon18] [regular polygon, regular polygon sides=8, draw] {};
\draw[thick] (Octagon78.corner 3) -- (Octagon78.corner 2);

\draw[->,shorten <=1mm,shorten >=1mm] (Octagon12) to (Octagon13);
\draw[->,shorten <= 1mm, shorten >= 1mm] (Octagon13) to (Octagon14);
\draw[->,shorten <= 1mm, shorten >= 1mm] (Octagon13) to (Octagon23);
\draw[->,shorten <= 1mm, shorten >= 1mm] (Octagon14) to (Octagon15);
\draw[->,shorten <= 1mm, shorten >= 1mm] (Octagon14) to (Octagon34);
\draw[->,shorten <= 1mm, shorten >= 1mm] (Octagon15) to (Octagon16);
\draw[->,shorten <= 1mm, shorten >= 1mm] (Octagon15) to (Octagon45);
\draw[->,shorten <= 1mm, shorten >= 1mm] (Octagon16) to (Octagon17);
\draw[->,shorten <= 1mm, shorten >= 1mm] (Octagon16) to (Octagon56);
\draw[->,shorten <= 1mm, shorten >= 1mm] (Octagon17) to (Octagon18);
\draw[->,shorten <= 1mm, shorten >= 1mm] (Octagon17) to (Octagon67);

\draw[->,shorten <= 1mm, shorten >= 1mm] (Octagon34) to (Octagon13);
\draw[->,shorten <= 1mm, shorten >= 1mm] (Octagon45) to (Octagon14);
\draw[->,shorten <= 1mm, shorten >= 1mm] (Octagon56) to (Octagon15);
\draw[->,shorten <= 1mm, shorten >= 1mm] (Octagon67) to (Octagon16);
\draw[->,shorten <= 1mm, shorten >= 1mm] (Octagon78) to (Octagon17);

\end{tikzpicture}}
\end{center}
We see that for any $n$ this quantum cluster algebra is of type $A_{n-3}$, since the subquiver on the vertices $\{\mathbf{13},\ldots,\mathbf{1(\itbf{n}-1)}\}$ is an orientation of the Dynkin diagram of this type.

The quasi-commutation matrix $L$ has four blocks, corresponding to the two blocks of $B$:
\begin{align*}
  L_{[1i][1k]} & =
     \begin{cases}
       -1 & \text{for}\ i\geq k+1 \\
        0 & \text{for}\ i=k \\
        1 & \text{for}\ i\leq k-1
     \end{cases} &
  L_{[1i][l(l+1)]} & =
     \begin{cases}
        0 & \text{for}\ i\geq l+2 \\
        1 & \text{for}\ i=l\ \text{or}\ l+1 \\
        2 & \text{for}\ i\leq l-1
     \end{cases} \\
  L_{[j(j+1)][1k]} & =
     \begin{cases}
        0 & \text{for}\ j\leq k-2 \\
       -1 & \text{for}\ j=k-1\ \text{or}\ k \\
       -2 & \text{for}\ j\geq k+1
     \end{cases} &
  L_{[j(j+1)][l(l+1)]} & =
     \begin{cases}
        2 & \text{for}\ j\leq l-2 \\
        1 & \text{for}\ j=l-1 \\
        0 & \text{for}\ j=l \\
       -1 & \text{for}\ j=l+1 \\
       -2 & \text{for}\ j\geq l+2
     \end{cases} 
\end{align*}
for $2\leq i,k\leq n$ and $1\leq j,l\leq n-1$.  These values may be verified easily, using the well-known quasi-commutation relations for quantum minors (see, for example, \cite{Scott-QMinors}).  For $n=8$, this matrix is
\[ L = 
\bordermatrix{
 & \scriptstyle{[18]} & \scriptstyle{[17]} & \scriptstyle{[16]} & \scriptstyle{[15]} & \scriptstyle{[14]} & \scriptstyle{[13]} & \scriptstyle{[12]} & \scriptstyle{[23]} & \scriptstyle{[34]} & \scriptstyle{[45]} & \scriptstyle{[56]} & \scriptstyle{[67]} & \scriptstyle{[78]} \cr
\scriptstyle{[18]} & 0 & -1 & -1 & -1 & -1 & -1 & -1 & 0 & 0 & 0 & 0 & 0 & 1 \cr
\scriptstyle{[17]} & 1 & \mathbf{0} & \mathbf{-1} & \mathbf{-1} & \mathbf{-1} & \mathbf{-1} & -1 & 0 & 0 & 0 & 0 & 1 & 1 \cr
\scriptstyle{[16]} & 1 & \mathbf{1} & \mathbf{0} & \mathbf{-1} & \mathbf{-1} & \mathbf{-1} & -1 & 0 & 0 & 0 & 1 & 1 & 2 \cr
\scriptstyle{[15]} & 1 & \mathbf{1} & \mathbf{1} & \mathbf{0} & \mathbf{-1} & \mathbf{-1} & -1 & 0 & 0 & 1 & 1 & 2 & 2 \cr
\scriptstyle{[14]} & 1 & \mathbf{1} & \mathbf{1} & \mathbf{1} & \mathbf{0} & \mathbf{-1} & -1 & 0 & 1 & 1 & 2 & 2 & 2 \cr
\scriptstyle{[13]} & 1 & \mathbf{1} & \mathbf{1} & \mathbf{1} & \mathbf{1} & \mathbf{0} & -1 & 1 & 1 & 2 & 2 & 2 & 2 \cr
\scriptstyle{[12]} & 1 & 1 & 1 & 1 & 1 & 1 &  0 &  1 &  2 &  2 &  2 &  2 & 2 \cr
\scriptstyle{[23]} & 0 & 0 & 0 & 0 & 0 & -1 & -1 &  0 &  1 &  2 &  2 &  2 & 2 \cr
\scriptstyle{[34]} & 0 & 0 & 0 & 0 & -1 & -1 & -2 & -1 &  0 &  1 &  2 &  2 & 2 \cr
\scriptstyle{[45]} &  0 &  0 &  0 & -1 & -1 & -2 & -2 & -2 & -1 &  0 &  1 &  2 & 2 \cr
\scriptstyle{[56]} &  0 &  0 & -1 & -1 & -2 & -2 & -2 & -2 & -2 & -1 &  0 &  1 & 2 \cr
\scriptstyle{[67]} &  0 & -1 & -1 & -2 & -2 & -2 & -2 & -2 & -2 & -2 & -1 &  0 & 1 \cr
\scriptstyle{[78]} & -1 & -1 & -2 & -2 & -2 & -2 & -2 & -2 & -2 & -2 & -2 & -1 & 0 
} \]
We do not give here a picture of the quiver $\Gamma(L)$, as it is rather unwieldy and does not give any additional information.

We claim that $B$ and $L$ are compatible.

\begin{proposition}\label{Gr-mats-compat} $B^{T}L=(0_{n-3,1}\ 2I_{n-3}\ 0_{n-3,n-1})$.
\end{proposition}

\begin{proof}
We separate the proof into two calculations, corresponding to the two blocks of $B$:

\begin{description}
\item[Block 1:] \textit{Claim:} $(B^{T}L)_{[1a][1b]}=2\delta_{ab}$ for $n-1\geq a\geq 3$ and $n\geq b\geq 2$.
\begin{eqnarray*}
    (B^{T}L)_{[1a][1b]} & = & \sum_{[1r]=[1n]}^{[12]} B_{[1r][1a]}L_{[1r][1b]}+\sum_{[r(r+1)]=[23]}^{[(n-1)n]} B_{[r(r+1)][1a]}L_{[r(r+1)][1b]} \\
      & \stackrel{(\star)}{=} & (B_{[1(a+1)][1a]}L_{[1(a+1)][1b]}+B_{[1(a-1)][1a]}L_{[1(a-1)][1b]}) \\
      & & \qquad +(B_{[(a-1)a][1a]}L_{[(a-1)a][1b]}+B_{[a(a+1)][1a]}L_{[a(a+1)][1b]}) \\
      & = & (-1)L_{[1(a+1)][1b]}+L_{[1(a-1)][1b]}+(-1)L_{[(a-1)a][1b]}+L_{[a(a+1)][1b]} \\
      & = & 
         \begin{cases}
          (-1)(-1)+(-1)+(-1)(-2)+(-2) & \text{if}\ a\geq b+2 \\
          (-1)(-1)+(0) +(-1)(-1)+(-2) & \text{if}\ a=b+1 \\
          (-1)(-1)+(1) +(-1)(-1)+(-1) & \text{if}\ a=b \\
          (-1)(0) +(1) +(-1)(0) +(-1) & \text{if}\ a=b-1 \\
          (-1)(1) +(1) +(-1)(0) +(0)  & \text{if}\ a\leq b-2
         \end{cases} \\
      & = & 
         \begin{cases}
          0 & \text{if}\ a\geq b+2 \\
          0 & \text{if}\ a=b+1 \\
          2 & \text{if}\ a=b \\
          0 & \text{if}\ a=b-1 \\
          0 & \text{if}\ a\leq b-2
         \end{cases} \\
      & = & 2\delta_{ab}
\end{eqnarray*}
where equality $(\star)$ holds since all other $B$-entries are zero.
\item[Block 2:] \textit{Claim:} $(B^{T}L)_{[1a][c(c+1)]}=0$ for $n-1\geq a\geq 3$ and $2\leq c \leq n-1$.
\begin{eqnarray*}
   (B^{T}L)_{[1a][c(c+1)]} & = & \sum_{[1r]=[1n]}^{[12]} B_{[1r][1a]}L_{[1r][c(c+1)]}+\sum_{[r(r+1)]=[23]}^{[(n-1)n]} B_{[r(r+1)][1a]}L_{[r(r+1)][c(c+1)]} \\
      & = & (-1)L_{[1(a+1)][c(c+1)]}+L_{[1(a-1)][c(c+1)]} \\ 
      & & \qquad +(-1)L_{[(a-1)a][c(c+1)]}+L_{[a(a+1)][c(c+1)]} \\
      & = & 
         \begin{cases}
          (-1)(0)+(0)+(-1)(-2)+(-2) & \text{if}\ a\geq c+3 \\
          (-1)(0)+(1)+(-1)(-1)+(-2) & \text{if}\ a=c+2 \\
          (-1)(0)+(1)+(-1)(0) +(-1) & \text{if}\ a=c+1 \\
          (-1)(1)+(2)+(-1)(1) +(0)  & \text{if}\ a=c \\
          (-1)(1)+(2)+(-1)(2) +(1)  & \text{if}\ a=c-1 \\
          (-1)(2)+(2)+(-1)(2) +(2)  & \text{if}\ a\leq c-2
         \end{cases} \\
      & = & 0 
\end{eqnarray*} 
\end{description}
Hence $B^{T}L$ has the stated form.
\end{proof}

In order to show that $\mathbb{K}_{q}[\mathrm{Gr}(2,n)]$ is a quantum cluster algebra, we must demonstrate that iterated mutation produces cluster variables that are in this algebra (rather than requiring any localisation) and that every generator occurs in some cluster.  We will do this inductively.

\begin{proposition} Let $3\leq \alpha<\beta \leq n-1$.  After the sequence of mutations\footnote{We index the mutations by the same labels as the elements they mutate, rather than by position in the cluster, for ease of identification of the result of mutation.  The compromise with this choice is that the label differs as more mutations are performed: the mutation label is that by which the element in that position was known before that mutation.  We will compose right to left, so that $\mu_{[1i]}\circ \mu_{[1(i+1)]}$ will mean perform $\mu_{[1(i+1)]}$ (replacing $\qminorbf{1(\itbf{i}+1)}$ by $\qminorbf{\itbf{i}(\itbf{i}+2)}$) then $\mu_{[1i]}$.  (So the $[1(i+1)]$-position in the cluster, to the left of $[1i]$, is in fact filled by $\qminorbf{\itbf{i}(\itbf{i}+2)}$ but we do not write $\mu_{[i(i+2)]}$ for $\mu_{[1(i+1)]}$ unless we mean to mutate $\qminorbf{\itbf{i}(\itbf{i}+2)}$.} $\mu_{[1(\alpha+1)]}\circ \mu_{[1(\alpha+2)]}\circ \cdots \circ \mu_{[1(\beta-1)]}$, the cluster obtained contains the mutable variables 
\[ (\qminorbf{1(\itbf{n}-1)},\ldots,\qminorbf{1\bbeta},\qminorbf{(\bbeta-2)\bbeta},\qminorbf{(\bbeta-3)\bbeta},\ldots,\qminorbf{\balpha \bbeta},\qminorbf{1\balpha},\qminorbf{1(\balpha-1)},\ldots,\qminorbf{13}). \]
Furthermore, the exchange matrix $B$ of this cluster has only the following non-zero entries in its $[1\alpha]$-column: $B_{[\alpha \beta][1\alpha]}=B_{[1(\alpha-1)][1\alpha]}=1$ and $B_{[1\beta][1\alpha]}=B_{[(\alpha-1)\alpha][1\alpha]}=-1$.
\end{proposition}

\begin{proof}
We work by induction on $\alpha$ and compute the mutation $\mu_{[1\alpha]}$ applied to the cluster and exchange matrix in the statement above.  This mutation exchanges $\qminorbf{1\balpha}$ for some new element $X$ which is determined by the quantum exchange relation.  This relation is determined by the non-zero entries of $B$ stated above and so is
\begin{eqnarray*} 
X & = & M\!\!\!\bordermatrix{ & \scriptscriptstyle{\cdots} & \scriptscriptstyle{[(\alpha+1)\beta]} & \scriptscriptstyle{[\alpha \beta]} & \scriptscriptstyle{[1\alpha]} & \scriptscriptstyle{[1(\alpha-1)]} & \scriptscriptstyle{[1(\alpha-2)]} & \scriptscriptstyle{\cdots} & \cr & \cdots & 0 & 1 & -1 & 1 & 0 & \cdots } \\ 
  &  & +M\!\!\!\bordermatrix{ & & & \scriptscriptstyle{[1\beta]} & & & & \scriptscriptstyle{[1\alpha]} & & & & \scriptscriptstyle{[(\alpha-1)\alpha]} & & \cr  & \,\cdots & 0 & 1 & 0 & \cdots & 0 & -1 & 0 & \cdots & 0 & 1 & 0 & \cdots } \\
  & = & q^{r}\qminorbf{\balpha \bbeta}(\qminorbf{1\balpha})^{-1}\qminorbf{1(\balpha-1)}+q^{s}\qminorbf{1\bbeta}(\qminorbf{1\balpha})^{-1}\qminorbf{(\balpha-1)\balpha}
\end{eqnarray*}
where $r$ and $s$ are integers to be calculated from $B$ and the quasi-commutation matrix.  The coefficient associated to the monomial $M(a_{1},\ldots,a_{r})$ is $\frac{1}{2}\sum_{i<j} a_{i}a_{j}L_{ji}$, with $L_{ji}$ the $(j,i)$-entry of the quasi-commutation matrix $L$.  Hence the integers $r$ and $s$ have the following values:
\begin{align*}
r & = \textstyle{\frac{1}{2}}(1\cdot (-1)\cdot L_{[1\alpha][\alpha \beta]}+1\cdot 1\cdot L_{[1(\alpha-1)][1\alpha]}+(-1)\cdot 1\cdot L_{[1(\alpha-1)][1\alpha]}) \\
 & = \textstyle{\frac{1}{2}}(1-2+1) \\
 & = 0 \\
\intertext{and}
s & = \textstyle{\frac{1}{2}}(1\cdot (-1)\cdot L_{[1\alpha][1\beta]}+1\cdot 1\cdot L_{[(\alpha-1)\alpha][1\beta]}+(-1)\cdot 1\cdot L_{[(\alpha-1)\alpha][1\alpha]}) \\
 & = \textstyle{\frac{1}{2}}(1+0-1) \\
 & = 0.
\end{align*}
Substituting these and applying quasi-commutation relations to move the $(\qminorbf{1\balpha})^{-1}$ terms to the left before re-arranging, we obtain the following equality involving $X$:
\[ \qminorbf{1\balpha}X = q^{-1}\qminorbf{1(\balpha-1)}\qminorbf{\balpha \bbeta}+q\qminorbf{1\bbeta}\qminorbf{(\balpha-1)\balpha}. \]
But the right-hand side of this equation is equal in $\mathbb{K}_{q}[\mathrm{Gr}(2,n)]$ to $\qminorbf{1\balpha}\qminorbf{(\balpha-1)\bbeta}$, by the quantum Pl\"{u}cker relations (see for instance \cite{KLR}).  Hence we deduce that $X=\qminorbf{(\balpha-1)\bbeta}$.

It remains to show that applying the mutation $\mu_{[1\alpha]}$ to $B$ gives a matrix $B'$ whose $[1(\alpha-1)]$-column has the correct entries.  To calculate this, the only information we need is that contained in the $[1\alpha]$-row from column $[1\alpha]$ onwards, the $[1\alpha]$-column of $B$ and the $[1(\alpha-1)]$-column of $B$.  The non-zero entries in this partial row and these columns are
\begin{enumerate}[label=(\alph*)]
\item $[1\alpha]$-row: $B_{[1\alpha][1(\alpha-1)]}=-1$,
\item $[1\alpha]$-column: $B_{[\alpha \beta][1\alpha]}=B_{[1(\alpha-1)][1\alpha]}=1$, $B_{[1\beta][1\alpha]}=B_{[(\alpha-1)\alpha][1\alpha]}=-1$,
\item $[1(\alpha-1)]$-column: \parbox[t]{18em}{$B_{[1(\alpha-2)][1(\alpha-1)]}=B_{[(\alpha-1)\alpha][1(\alpha-1)]}=1$, \hfill \\ $B_{[1\alpha][1(\alpha-1)]}=B_{[(\alpha-2)(\alpha-1)][1(\alpha-1)]}=-1$.\rule{0em}{1.2em}}
\end{enumerate}
The first and last of these are the same as for the initial exchange matrix, as our mutation sequence has only affected the exchange matrix at and to the left of the $[1\alpha]$-column.  The middle of these is exactly the claim on the form of the exchange matrix in the inductive hypothesis.

Recall that the formula for the matrix mutation $\mu_{k}$ applied to a matrix $B$ is 
\[ (\mu_{k}(B))_{ij} = \begin{cases} -B_{ij} & \text{if}\ i=k\ \text{or}\ j=k \\ B_{ij}+\frac{|B_{ik}|B_{kj}+B_{ik}|B_{kj}|}{2} & \text{otherwise} \end{cases} \]
As a result, in $B'$ entries in $[1(\alpha-1)]$-column will certainly not change whenever the $[1\alpha]$-column in $B$ has a zero, except for in the $[1\alpha]$-row, whose sign is changed.  (If $B_{ik}$ and $B_{kj}$ have equal magnitude but opposite signs then $B_{ij}$ will also not change.)  Thus it suffices to make the following set of calculations:
\begin{align}
 B'_{[1\beta][1(\alpha-1)]} & = B_{[1\beta][1(\alpha-1)]}+\frac{|B_{[1\beta][1\alpha]}|B_{[1\alpha][1(\alpha-1)]}+B_{[1\beta][1\alpha]}|B_{[1\alpha][1(\alpha-1)]}|}{2} \\ & = 0+\frac{|-1|(-1)+(-1)|-1|}{2} \notag \\ & = -1 \notag \\
 B'_{[\alpha \beta][1(\alpha-1)]} & = 0+\frac{|1|(-1)+(1)|-1|}{2} \\ & = 0 \notag \\
 B'_{[(\alpha-1)\beta][1(\alpha-1)]} & = -B_{[1\alpha][1(\alpha-1)]} \\ & = 1 \notag \\
\intertext{(This row in $B'$, the $[(\alpha-1)\beta]$-row, replaces the $[1\alpha]$-row in $B$; it is the latter's negative, by the mutation rule.)}
 B'_{[1(\alpha-1)][1(\alpha-1)]} & = 0+\frac{|1|(-1)+(1)|-1|}{2} \\ & = 0 \notag \\
 B'_{[1(\alpha-2)][1(\alpha-1)]} & = 1+\frac{|0|(-1)+(0)|-1|}{2} \\ & = 1 \notag \\
\intertext{(This entry and the one following are unchanged in $B'$ because there was a 0 in the adjacent $[1\alpha]$-column, but are included here explicitly for completeness, as they are non-zero entries in the $[1(\alpha-1)]$-column on $B'$.)}
 B'_{[(\alpha-2)(\alpha-1)][1(\alpha-1)]} & = -1+\frac{|0|(-1)+(0)|-1|}{2} \\ & = -1 \notag \\
 B'_{[(\alpha-1)(\alpha)][1(\alpha-1)]} & = 1+\frac{|-1|(-1)+(-1)|-1|}{2} \\ & = 0 \notag
\end{align}
From these, we see that the $[1(\alpha-1)]$-column of $B'$ has only four non-zero entries and these are $B'_{[(\alpha-1) \beta][1(\alpha-1)]}=B'_{[1(\alpha-2)][1(\alpha-1)]}=1$ and $B'_{[1\beta][1(\alpha-1)]}=B'_{[(\alpha-2)(\alpha-1)][1(\alpha-1)]}=-1$.\rule[-0.8em]{0em}{1.2em}  These are exactly as claimed for $B$ in the statement of this proposition, except with $\alpha-1$ instead of $\alpha$, \ie the inductive hypothesis for the exchange matrix has been verified.
\end{proof}

We can now complete the proof of our claim.

\begin{theorem}
\label{theo-grass} The quantum Grassmannian $\mathbb{K}_{q}[\mathrm{Gr}(2,n)]$ is a quantum cluster algebra of type $A_{n-3}$.
\end{theorem}

\begin{proof}  It follows immediately from the previous proposition that every generator of $\mathbb{K}_{q}[\mathrm{Gr}(2,n)]$ occurs in some cluster and can be obtained by a finite sequence of mutations from our given initial seed.  Furthermore, the general theory of (quantum) cluster algebras of type $A$ tells us that a cluster algebra of type $A_{l}$ has precisely $\frac{l(l+1)}{2}+l=\frac{l^{2}+3l}{2}$ mutable cluster variables (as this is the number of almost positive roots in this type).  For $l=n-3$, this number is $\frac{(n-3)(n-2)}{2}+(n-3)=\frac{n^{2}-3n}{2}$.  Hence since $|\curly{P}_{q}|-n=\frac{(n-1)n}{2}-n=\frac{n^{2}-3n}{2}$ also, the set of all mutable cluster variables must be equal to the set $\curly{P}_{q}\setminus \{\qminor{12},\ldots,\qminor{(n-1)n},\qminor{1n}\}$, since the latter is contained in the former by the proposition.  

In other words, taking all cluster variables (mutable ones and coefficients) together gives us exactly the set $\curly{P}_{q}$ of quantum Pl\"{u}cker coordinates.  From this we deduce that every cluster variable is a genuine quantum minor (\ie no localisation is required) and hence the subalgebra generated by all cluster variables is exactly $\mathbb{K}_{q}[\mathrm{Gr}(2,n)]$, which is therefore a quantum cluster algebra.
\end{proof}

\begin{remark} This proof uses in a critical way the fact that the quantum cluster algebra structure we have found is of finite type, and this argument will not extend to infinite types.
\end{remark}

\begin{remark} 
\label{qhypothesis} Although we have assumed that $q^{1/2}$ exists in $\mathbb{K}$, it turns out that this assumption is not needed. Indeed, it follows from the proof of Theorem \ref{theo-grass} that the powers of $q$ that appear in the (quantum) mutation are all integral numbers.
\end{remark}

We reproduce in Figure~\ref{Gr25-ex-graph} the diagram in \cite{Fomin-Reading} showing the exchange graph in the case $n=5$, with clusters identified with triangulations of a pentagon in the manner described previously.  The top vertex corresponds to the initial cluster described here, whose mutable variables are the quantum minors $\qminorbf{14}$ and $\qminorbf{13}$: we number the pentagon's vertices starting with 1 at the top and increasing clockwise.

\begin{figure}[ht]
\begin{center} 
\scalebox{0.75}{\begin{tikzpicture}
\node (Pentagon1) at (90:4cm) [regular polygon, regular polygon sides=5, draw,minimum size=2cm] {};
\node (Pentagon2) at (72+90:4cm) [regular polygon, regular polygon sides=5, draw,minimum size=2cm] {};
\node (Pentagon3) at (144+90:4cm) [regular polygon, regular polygon sides=5, draw,minimum size=2cm] {};
\node (Pentagon4) at (216+90:4cm) [regular polygon, regular polygon sides=5, draw,minimum size=2cm] {};
\node (Pentagon5) at (288+90:4cm) [regular polygon, regular polygon sides=5, draw,minimum size=2cm] {};
\draw[shorten <=3mm,shorten >=3mm] (Pentagon1) -- (Pentagon2);
\draw[shorten <=3mm,shorten >=3mm] (Pentagon2) -- (Pentagon3);
\draw[shorten <=3mm,shorten >=3mm] (Pentagon3) -- (Pentagon4);
\draw[shorten <=3mm,shorten >=3mm] (Pentagon4) -- (Pentagon5);
\draw[shorten <=3mm,shorten >=3mm] (Pentagon5) -- (Pentagon1);
\draw[] (Pentagon1.corner 1) -- (Pentagon1.corner 4);
\draw[] (Pentagon1.corner 1) -- (Pentagon1.corner 3);
\draw[] (Pentagon2.corner 1) -- (Pentagon2.corner 4);
\draw[] (Pentagon2.corner 2) -- (Pentagon2.corner 4);
\draw[] (Pentagon3.corner 2) -- (Pentagon3.corner 4);
\draw[] (Pentagon3.corner 2) -- (Pentagon3.corner 5);
\draw[] (Pentagon4.corner 2) -- (Pentagon4.corner 5);
\draw[] (Pentagon4.corner 3) -- (Pentagon4.corner 5);
\draw[] (Pentagon5.corner 3) -- (Pentagon5.corner 5);
\draw[] (Pentagon5.corner 1) -- (Pentagon5.corner 3);

\end{tikzpicture}}
\end{center}
\caption{Exchange graph for cluster algebra structure on $\complex[\mathrm{Gr}(2,5)]$ and its quantum analogue.}\label{Gr25-ex-graph}
\end{figure}

\subsection{The quantum Grassmannians \texorpdfstring{$\mathbb{K}_{q}[\mathrm{Gr}(3,n)]$, $n=6,7,8$}{Kq[Gr(3,n)], n=6,7,8}}\label{ss:Gr3-678}

We now consider the remaining quantum Grassmannians that should have finite type as quantum cluster algebras, namely $\mathbb{K}_{q}[\mathrm{Gr}(3,n)]$ for $n=6,7,8$.  For these examples, we give an initial quantum seed and a list of the quantum cluster variables obtained by certain mutations of this.  This list will be a list of all the quantum cluster variables (some of which are not quantum minors, as happens classically) and will be the quantum analogue of the corresponding table among Tables~1-4 in \cite{Scott-Grassmannians}.

\subsubsection{$\mathbb{K}_{q}[\mathrm{Gr}(3,6)]$}

For our initial quantum cluster we choose
\[ \tilde{\underline{y}}=(\qminorbf{145},\qminorbf{134},\qminorbf{125},\qminorbf{124},\qminor{123},\qminor{234},\qminor{345},\qminor{456},\qminor{156},\qminor{126}). \]  As above, the coefficients correspond to $k$-tuples of adjacent vertices of an $n$-gon---i.e. tuples of vertices of a hexagon, wrapping round and taking $\qminor{156}$ for $\{5,6,1\}$ and so on.  The exchange matrix for this cluster is most compactly described by its associated quiver.  This is as follows, with minor labels and with hexagons.

\vspace{1em}
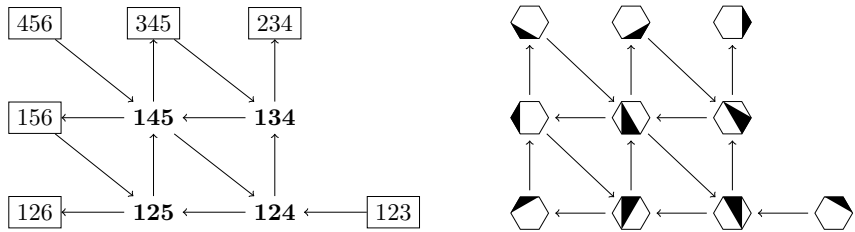
\begin{center}
\scalebox{0.825}{\begin{tikzpicture}
\node (145) at (0,0) {$\mathbf{145}$};
\node (125) [below=of 145] {$\mathbf{125}$};
\node (134) [right=of 145] {$\mathbf{134}$};
\node (124) [below=of 134] {$\mathbf{124}$};
\node (156) [left=of 145,rectangle,draw=black] {156};
\node (126) [below=of 156,rectangle,draw=black] {126};
\node (123) [right=of 124,rectangle,draw=black]  {123};
\node (456) [above=of 156,rectangle,draw=black] {456};
\node (345) [above=of 145,rectangle,draw=black] {345};
\node (234) [above=of 134,rectangle,draw=black] {234};

\draw[->] (124) to (125);
\draw[->] (125) to (126);

\draw[->] (134) to (145);
\draw[->] (145) to (156);

\draw[->] (124) to (134);
\draw[->] (134) to (234);

\draw[->] (125) to (145);
\draw[->] (145) to (345);

\draw[->] (156) to (125);
\draw[->] (456) to (145);

\draw[->] (145) to (124);
\draw[->] (345) to (134);

\draw[->] (123) to (124);

\node (GHexagon156) at (6,0) [regular polygon, regular polygon sides=6, draw,minimum size=0.6cm] {};
\fill (GHexagon156.corner 2) -- (GHexagon156.corner 4) -- (GHexagon156.corner 3);

\node (GHexagon126) [below=of GHexagon156,regular polygon, regular polygon sides=6, draw,minimum size=0.6cm] {};
\fill (GHexagon126.corner 2) -- (GHexagon126.corner 1) -- (GHexagon126.corner 3);

\node (GHexagon145) [right=of GHexagon156,regular polygon, regular polygon sides=6, draw,minimum size=0.6cm] {};
\fill (GHexagon145.corner 2) -- (GHexagon145.corner 5) -- (GHexagon145.corner 4);

\node (GHexagon125) [below=of GHexagon145,regular polygon, regular polygon sides=6, draw,minimum size=0.6cm] {};
\fill (GHexagon125.corner 2) -- (GHexagon125.corner 1) -- (GHexagon125.corner 4);

\node (GHexagon134) [right=of GHexagon145,regular polygon, regular polygon sides=6, draw,minimum size=0.6cm] {};
\fill (GHexagon134.corner 2) -- (GHexagon134.corner 6) -- (GHexagon134.corner 5);

\node (GHexagon124) [below=of GHexagon134,regular polygon, regular polygon sides=6, draw,minimum size=0.6cm] {};
\fill (GHexagon124.corner 2) -- (GHexagon124.corner 1) -- (GHexagon124.corner 5);

\node (GHexagon123) [right=of GHexagon124,regular polygon, regular polygon sides=6, draw,minimum size=0.6cm]  {};
\fill (GHexagon123.corner 2) -- (GHexagon123.corner 1) -- (GHexagon123.corner 6);

\node (GHexagon456) [above=of GHexagon156,regular polygon, regular polygon sides=6, draw,minimum size=0.6cm] {};
\fill (GHexagon456.corner 5) -- (GHexagon456.corner 4) -- (GHexagon456.corner 3);

\node (GHexagon345) [above=of GHexagon145,regular polygon, regular polygon sides=6, draw,minimum size=0.6cm] {};
\fill (GHexagon345.corner 6) -- (GHexagon345.corner 5) -- (GHexagon345.corner 4);

\node (GHexagon234) [above=of GHexagon134,regular polygon, regular polygon sides=6, draw,minimum size=0.6cm] {};
\fill (GHexagon234.corner 1) -- (GHexagon234.corner 6) -- (GHexagon234.corner 5);

\draw[->,shorten <= 1mm, shorten >= 1mm] (GHexagon124) to (GHexagon125);
\draw[->,shorten <= 1mm, shorten >= 1mm] (GHexagon125) to (GHexagon126);

\draw[->,shorten <= 1mm, shorten >= 1mm] (GHexagon134) to (GHexagon145);
\draw[->,shorten <= 1mm, shorten >= 1mm] (GHexagon145) to (GHexagon156);

\draw[->,shorten <= 1mm, shorten >= 1mm] (GHexagon124) to (GHexagon134);
\draw[->,shorten <= 1mm, shorten >= 1mm] (GHexagon134) to (GHexagon234);

\draw[->,shorten <= 1mm, shorten >= 1mm] (GHexagon125) to (GHexagon145);
\draw[->,shorten <= 1mm, shorten >= 1mm] (GHexagon145) to (GHexagon345);

\draw[->,shorten <= 1mm, shorten >= 1mm] (GHexagon126) to (GHexagon156);
\draw[->,shorten <= 1mm, shorten >= 1mm] (GHexagon156) to (GHexagon456);

\draw[->,shorten <= 1mm, shorten >= 1mm] (GHexagon156) to (GHexagon125);
\draw[->,shorten <= 1mm, shorten >= 1mm] (GHexagon456) to (GHexagon145);

\draw[->,shorten <= 1mm, shorten >= 1mm] (GHexagon145) to (GHexagon124);
\draw[->,shorten <= 1mm, shorten >= 1mm] (GHexagon345) to (GHexagon134);

\draw[->,shorten <= 1mm, shorten >= 1mm] (GHexagon123) to (GHexagon124);

\end{tikzpicture}}
\end{center}
This quiver is mutation-equivalent to the Dynkin diagram of type $D_{4}$. 
The quasi-commutation matrix associated to $\tilde{\underline{y}}$ is
\[ L = 
\bordermatrix{
 & \scriptstyle{[145]} & \scriptstyle{[134]} & \scriptstyle{[125]} & \scriptstyle{[124]} & \scriptstyle{[123]} & \scriptstyle{[234]} & \scriptstyle{[345]} & \scriptstyle{[456]} & \scriptstyle{[156]} & \scriptstyle{[126]} \cr
\scriptstyle{[145]} & \mathbf{0} & \mathbf{-1} & \mathbf{-1} & \mathbf{-1} & -2 & 0 & 1 & 1 & 1 & 0 \cr
\scriptstyle{[134]} & \mathbf{1} & \mathbf{0} & \mathbf{0} & \mathbf{-1} & -1 & 1 & 1 & 2 & 2 & 0 \cr
\scriptstyle{[125]} & \mathbf{1} & \mathbf{0} & \mathbf{0} & \mathbf{-1} & -1 & 0 & 2 & 2 & 1 & 1 \cr
\scriptstyle{[136]} & \mathbf{1} & \mathbf{1} & \mathbf{1} & \mathbf{0} & -1 & 1 & 2 & 2 & 2 & 1 \cr
\scriptstyle{[124]} & 2 & 1 & 1 & 1 & 0 & 1 & 2 & 3 & 2 & 1 \cr
\scriptstyle{[234]} & 0 & -1 & 0 & -1 & -1 & 0 & 1 & 2 & 1 & 0 \cr
\scriptstyle{[345]} & -1 & -1 & -2 & -2 & -2 & -1 &  0 &  1 &  0 &  -1 \cr
\scriptstyle{[456]} & -1 & -2 & -2 & -2 & -3 & -2 & -1 & 0 &  -1 &  -2 \cr
\scriptstyle{[156]} & -1 & -2 & -1 & -2 & -2 & -1 & 0 & 1 &  0 &  -1 \cr
\scriptstyle{[126]} & 0 & 0 & -1 & -1 & -1 & 0 & 1 & 2 & 1 &  0 \cr
}. \]
These quasi-commutation relations may be derived from the formula in \cite{Scott-QMinors} and it is straightforward to verify that $B^{T}L$ has the block form $( 2I_{4}\ 0_{4,6})$ and hence $B$ and $L$ are compatible.

The complete list of mutable quantum cluster variables is then obtained by repeated mutation, to give the list in Table~\ref{Gr36}.  This was obtained with assistance from the computer program Magma (\cite{Magma}).

\renewcommand{\arraystretch}{1.25}
\begin{table}
\begin{center}
\begin{tabular}{|l|l|l|} \hline
Mutation & Quantum cluster variable obtained & Almost-positive root \\ \hline
 & $\qminor{145}$ & $[-1,0,0,0]$ \\ \hline
 & $\qminor{134}$ & $[0,1,0,1]$ \\ \hline
 & $\qminor{125}$ & $[0,1,1,0]$ \\ \hline
 & $\qminor{124}$ & $[0,1,1,1]$ \\ \hline
$\mu_{1}$ & $X_{q}^{123456} \defeq q^{3/2}(\qminor{356}\qminor{124}-q\qminor{456}\qminor{123})$ & $[1,2,1,1]$ \\ \hline
$\mu_{2}$ & $\qminor{245}$ & $[0,0,1,0]$ \\ \hline
$\mu_{3}$ & $\qminor{146}$ & $[0,0,0,1]$ \\ \hline
$\mu_{4}$ & $\qminor{135}$ & $[0,1,0,0]$ \\ \hline
$\mu_{2} \circ \mu_{1}$ & $\qminor{256}$ & $[1,1,1,0]$ \\ \hline
$\mu_{3} \circ \mu_{1}$ & $\qminor{346}$ & $[1,1,0,1]$ \\ \hline
$\mu_{4} \circ \mu_{1}$ & $\qminor{356}$ & $[1,1,0,0]$ \\ \hline
$\mu_{4} \circ \mu_{2}$ & $\qminor{235}$ & $[0,0,0,-1]$ \\ \hline
$\mu_{4} \circ \mu_{3}$ & $\qminor{136}$ & $[0,0,-1,0]$ \\ \hline
$\mu_{2} \circ \mu_{1} \circ \mu_{3}$ & $\qminor{246}$ & $[1,1,1,1]$ \\ \hline
$\mu_{4} \circ \mu_{3} \circ \mu_{2}$ &  $Y_{q}^{123456} \defeq q^{1/2}(\qminor{236}\qminor{145}-q^{2}\qminor{456}\qminor{123})$ & $[0,-1,0,0]$ \\ \hline
$\mu_{4} \circ \mu_{2} \circ \mu_{1} \circ \mu_{3}$ & $\qminor{236}$ & $[1,0,0,0]$ \\ \hline
\end{tabular}
\caption{Mutable quantum cluster variables for $\mathbb{K}_{q}[\mathrm{Gr}(3,6)]$}\label{Gr36}
\end{center}
\end{table}
\renewcommand{\arraystretch}{1}

We see that Table~\ref{Gr36} contains 16 mutable quantum cluster variables, 14 of which are quantum minors which taken together with the 6 coefficients yield the whole set of quantum Pl\"{u}cker coordinates generating $\mathbb{K}_{q}[\mathrm{Gr}(3,6)]$.  Hence $\mathbb{K}_{q}[\mathrm{Gr}(3,6)]$ is a quantum cluster algebra of type $D_{4}$.

\begin{remark} The elements $X_{q}^{123456}$ and $Y_{q}^{123456}$ are quantum analogues of the quadratic regular functions described by Scott in \cite[Theorem~6]{Scott-Grassmannians}.  There $X^{123456} \defeq \Delta^{134}\Delta^{256}-\Delta^{156}\Delta^{234}$, whereas we have set $X_{q}^{123456}=q^{3/2}(\qminor{356}\qminor{124}-q\qminor{456}\qminor{123})$.  In fact, \[ q^{-1/2}\qminor{134}\qminor{256}-q^{3/2}\qminor{156}\qminor{234}=q^{3/2}\qminor{356}\qminor{124}-q^{5/2}\qminor{456}\qminor{123}=X_{q}^{123456}; \] there will typically be several different choices of expression for each quantum cluster variable.
\end{remark}

\begin{remark} The almost-positive root bijection is not obtained directly from our choice of initial quantum cluster.  One must in fact choose a quantum cluster whose associated exchange quiver is the Dynkin graph determining the type with an orientation such that every mutable vertex is either a source or a sink.  Then according to the general theory of (quantum) cluster algebras, the bijection arises by associating negative simple roots to the initial mutable cluster variables and the positive roots correspond to variables obtained by mutation in the natural way.  We have used the same bijection as Scott (\cite{Scott-Grassmannians}), which arises from the cluster having exchange quiver as follows.

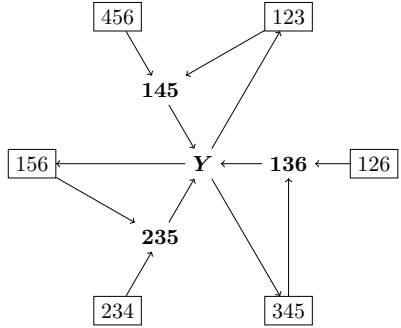
\begin{center}
\scalebox{0.75}{\begin{tikzpicture}
\node (Y) at (0:0) {$\mathitbf{Y}$};
\node (136) at (0:1.5cm) {\textbf{136}};
\node (145) at (120:1.5cm) {\textbf{145}};
\node (235) at (240:1.5cm) {\textbf{235}};

\node (123) [rectangle,draw=black] at (60:3cm) {123};
\node (126) [rectangle,draw=black] at (0:3cm) {126};
\node (156) [rectangle,draw=black] at (180:3cm) {156};
\node (456) [rectangle,draw=black] at (120:3cm) {456};
\node (345) [rectangle,draw=black] at (300:3cm) {345};
\node (234) [rectangle,draw=black] at (240:3cm) {234};

\draw[->] (136) to (Y);
\draw[->] (145) to (Y);
\draw[->] (235) to (Y);

\draw[->] (Y) to (123);
\draw[->] (123) to (145);

\draw[->] (234) to (235);

\draw[->] (Y) to (345);
\draw[->] (345) to (136);

\draw[->] (456) to (145);

\draw[->] (Y) to (156);
\draw[->] (156) to (235);

\draw[->] (126) to (136);

\end{tikzpicture}}
\end{center}
\end{remark}

\begin{remark} We note that in contrast to the situation for $\mathbb{K}_{q}[\mathrm{Gr}(2,n)]$, we cannot drop the requirement for there to be a square root of $q$ in $\mathbb{K}$, as it appears in the definitions of the cluster variables $Y_{q}^{123456}$ and $X_{q}^{123456}$.
\end{remark}

\subsubsection{$\mathbb{K}_{q}[\mathrm{Gr}(3,7)]$}

By analogy with the previous example, we choose for our initial cluster
\[ \tilde{\underline{y}}=(\qminorbf{156},\qminorbf{145},\qminorbf{134},\qminorbf{126},\qminorbf{125},\qminorbf{124},\qminor{123},\qminor{234},\qminor{345},\qminor{456},\qminor{567},\qminor{167},\qminor{127}). \]  The coefficients correspond to triples of adjacent vertices of a heptagon and the exchange matrix is as determined its associated quiver, which is as follows.
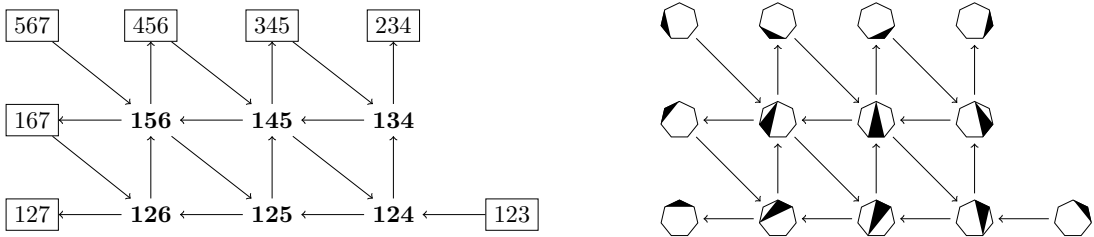
\begin{center}
\scalebox{0.825}{\begin{tikzpicture}
\node (156) at (0,0) {$\mathbf{156}$};
\node (126) [below=of 156] {$\mathbf{126}$};
\node (145) [right=of 156] {$\mathbf{145}$};
\node (125) [below=of 145] {$\mathbf{125}$};
\node (134) [right=of 145] {$\mathbf{134}$};
\node (124) [below=of 134] {$\mathbf{124}$};
\node (123) [right=of 124,rectangle,draw=black]  {123};
\node (127) [left=of 126,rectangle,draw=black] {127};
\node (167) [left=of 156,rectangle,draw=black] {167};
\node (567) [above=of 167,rectangle,draw=black] {567};
\node (456) [above=of 156,rectangle,draw=black] {456};
\node (345) [above=of 145,rectangle,draw=black] {345};
\node (234) [above=of 134,rectangle,draw=black] {234};

\draw[->] (124) to (125);
\draw[->] (125) to (126);
\draw[->] (126) to (127);

\draw[->] (134) to (145);
\draw[->] (145) to (156);
\draw[->] (156) to (167);

\draw[->] (124) to (134);
\draw[->] (134) to (234);

\draw[->] (125) to (145);
\draw[->] (145) to (345);

\draw[->] (126) to (156);
\draw[->] (156) to (456);

\draw[->] (167) to (126);
\draw[->] (567) to (156);

\draw[->] (156) to (125);
\draw[->] (456) to (145);

\draw[->] (145) to (124);
\draw[->] (345) to (134);

\draw[->] (123) to (124);


\node (GHeptagon156) at (10,0) [regular polygon, regular polygon sides=7, draw,minimum size=0.6cm] {};
\fill (GHeptagon156.corner 1) -- (GHeptagon156.corner 4) -- (GHeptagon156.corner 3);

\node (GHeptagon126) [below=of GHeptagon156,regular polygon, regular polygon sides=7, draw,minimum size=0.6cm] {};
\fill (GHeptagon126.corner 1) -- (GHeptagon126.corner 7) -- (GHeptagon126.corner 3);

\node (GHeptagon145) [right=of GHeptagon156,regular polygon, regular polygon sides=7, draw,minimum size=0.6cm] {};
\fill (GHeptagon145.corner 1) -- (GHeptagon145.corner 5) -- (GHeptagon145.corner 4);

\node (GHeptagon125) [below=of GHeptagon145,regular polygon, regular polygon sides=7, draw,minimum size=0.6cm] {};
\fill (GHeptagon125.corner 1) -- (GHeptagon125.corner 7) -- (GHeptagon125.corner 4);

\node (GHeptagon134) [right=of GHeptagon145,regular polygon, regular polygon sides=7, draw,minimum size=0.6cm] {};
\fill (GHeptagon134.corner 1) -- (GHeptagon134.corner 6) -- (GHeptagon134.corner 5);

\node (GHeptagon124) [below=of GHeptagon134,regular polygon, regular polygon sides=7, draw,minimum size=0.6cm] {};
\fill (GHeptagon124.corner 1) -- (GHeptagon124.corner 7) -- (GHeptagon124.corner 5);

\node (GHeptagon123) [right=of GHeptagon124,regular polygon, regular polygon sides=7, draw,minimum size=0.6cm]  {};
\fill (GHeptagon123.corner 1) -- (GHeptagon123.corner 7) -- (GHeptagon123.corner 6);

\node (GHeptagon127) [left=of GHeptagon126,regular polygon, regular polygon sides=7, draw,minimum size=0.6cm] {};
\fill (GHeptagon127.corner 1) -- (GHeptagon127.corner 7) -- (GHeptagon127.corner 2);

\node (GHeptagon167) [left=of GHeptagon156,regular polygon, regular polygon sides=7, draw,minimum size=0.6cm] {};
\fill (GHeptagon167.corner 1) -- (GHeptagon167.corner 3) -- (GHeptagon167.corner 2);

\node (GHeptagon567) [above=of GHeptagon167,regular polygon, regular polygon sides=7, draw,minimum size=0.6cm] {};
\fill (GHeptagon567.corner 4) -- (GHeptagon567.corner 3) -- (GHeptagon567.corner 2);

\node (GHeptagon456) [above=of GHeptagon156,regular polygon, regular polygon sides=7, draw,minimum size=0.6cm] {};
\fill (GHeptagon456.corner 5) -- (GHeptagon456.corner 4) -- (GHeptagon456.corner 3);

\node (GHeptagon345) [above=of GHeptagon145,regular polygon, regular polygon sides=7, draw,minimum size=0.6cm] {};
\fill (GHeptagon345.corner 6) -- (GHeptagon345.corner 5) -- (GHeptagon345.corner 4);

\node (GHeptagon234) [above=of GHeptagon134,regular polygon, regular polygon sides=7, draw,minimum size=0.6cm] {};
\fill (GHeptagon234.corner 7) -- (GHeptagon234.corner 6) -- (GHeptagon234.corner 5);

\draw[->,shorten <= 1mm, shorten >= 1mm] (GHeptagon124) to (GHeptagon125);
\draw[->,shorten <= 1mm, shorten >= 1mm] (GHeptagon125) to (GHeptagon126);
\draw[->,shorten <= 1mm, shorten >= 1mm] (GHeptagon126) to (GHeptagon127);

\draw[->,shorten <= 1mm, shorten >= 1mm] (GHeptagon134) to (GHeptagon145);
\draw[->,shorten <= 1mm, shorten >= 1mm] (GHeptagon145) to (GHeptagon156);
\draw[->,shorten <= 1mm, shorten >= 1mm] (GHeptagon156) to (GHeptagon167);

\draw[->,shorten <= 1mm, shorten >= 1mm] (GHeptagon124) to (GHeptagon134);
\draw[->,shorten <= 1mm, shorten >= 1mm] (GHeptagon134) to (GHeptagon234);

\draw[->,shorten <= 1mm, shorten >= 1mm] (GHeptagon125) to (GHeptagon145);
\draw[->,shorten <= 1mm, shorten >= 1mm] (GHeptagon145) to (GHeptagon345);

\draw[->,shorten <= 1mm, shorten >= 1mm] (GHeptagon126) to (GHeptagon156);
\draw[->,shorten <= 1mm, shorten >= 1mm] (GHeptagon156) to (GHeptagon456);

\draw[->,shorten <= 1mm, shorten >= 1mm] (GHeptagon167) to (GHeptagon126);
\draw[->,shorten <= 1mm, shorten >= 1mm] (GHeptagon567) to (GHeptagon156);

\draw[->,shorten <= 1mm, shorten >= 1mm] (GHeptagon156) to (GHeptagon125);
\draw[->,shorten <= 1mm, shorten >= 1mm] (GHeptagon456) to (GHeptagon145);

\draw[->,shorten <= 1mm, shorten >= 1mm] (GHeptagon145) to (GHeptagon124);
\draw[->,shorten <= 1mm, shorten >= 1mm] (GHeptagon345) to (GHeptagon134);

\draw[->,shorten <= 1mm, shorten >= 1mm] (GHeptagon123) to (GHeptagon124);

\end{tikzpicture}}
\end{center}
This quiver is mutation-equivalent to the Dynkin diagram of type $E_{6}$.  We omit the corresponding quasi-commutation matrix, which can be recovered easily, and simply note that the compatibility calculation gives the matrix $(2I_{6}\ 0_{6,7})$.

We will also not give the full list of quantum cluster variables.  By computer-aided calculation, we have verified that the set of quantum cluster variables in this example consists of all quantum Pl\"{u}cker coordinates together with the following 14 additional elements.  Let $\{a,b,c,d,e,f\}$ be a subset of $\{1,\ldots ,7 \}$ written in increasing order.  Define
\begin{align*}
X_{q}^{abcdef} & \defeq q^{3/2}(\qminor{cef}\qminor{abd}-q\qminor{def}\qminor{abc}) \\
\intertext{and}
Y_{q}^{abcdef} & \defeq q^{1/2}(\qminor{bcf}\qminor{ade}-q^{2}\qminor{def}\qminor{abc}).
\end{align*}
Again, these are quantizations of the quadratic regular functions described by Scott in the classical case.  (Also, if one takes $\{a,b,c,d,e,f\}=\{1,\ldots ,6\}$ one obtains the functions described for $n=6$ above.)

We again conclude that provided $q$ has a square root in $\mathbb{K}$, the quantum Grassmannian $\mathbb{K}_{q}[\mathrm{Gr}(3,7)]$ is a quantum cluster algebra of type $E_{6}$.

\subsubsection{$\mathbb{K}_{q}[\mathrm{Gr}(3,8)]$}

The choice of initial cluster in the same family as those above should now be clear.  That is, we take
\begin{multline*} \tilde{\underline{y}}=(\qminorbf{167},\qminorbf{156},\qminorbf{145},\qminorbf{134},\qminorbf{127},\qminorbf{126},\qminorbf{125},\qminorbf{124}, \\ \qminor{123},\qminor{234},\qminor{345},\qminor{456},\qminor{567},\qminor{678},\qminor{178},\qminor{128}).\end{multline*}
The coefficients correspond to triples of adjacent vertices of an octagon and the exchange matrix is as determined its associated quiver, which is as follows.
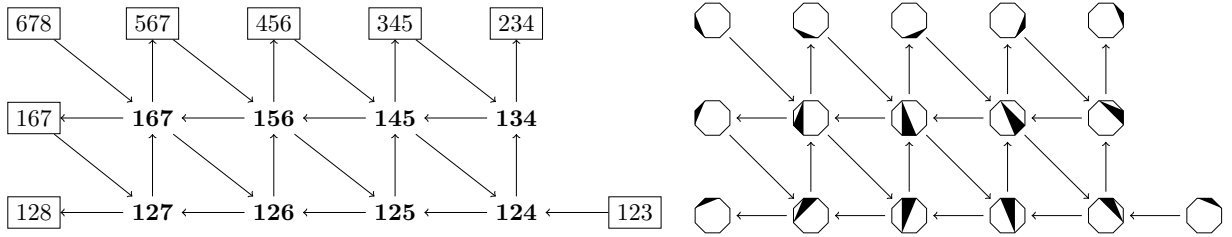
\begin{center}
\scalebox{0.825}{\begin{tikzpicture}
\node (167) at (0,0) {$\mathbf{167}$};
\node (127) [below=of 167] {$\mathbf{127}$};
\node (156) [right=of 167] {$\mathbf{156}$};
\node (126) [below=of 156] {$\mathbf{126}$};
\node (145) [right=of 156] {$\mathbf{145}$};
\node (125) [below=of 145] {$\mathbf{125}$};
\node (134) [right=of 145] {$\mathbf{134}$};
\node (124) [below=of 134] {$\mathbf{124}$};
\node (123) [right=of 124,rectangle,draw=black]  {123};
\node (128) [left=of 127,rectangle,draw=black] {128};
\node (178) [left=of 167,rectangle,draw=black] {167};
\node (678) [above=of 178,rectangle,draw=black] {678};
\node (567) [above=of 167,rectangle,draw=black] {567};
\node (456) [above=of 156,rectangle,draw=black] {456};
\node (345) [above=of 145,rectangle,draw=black] {345};
\node (234) [above=of 134,rectangle,draw=black] {234};

\draw[->] (123) to (124);
\draw[->] (124) to (125);
\draw[->] (125) to (126);
\draw[->] (126) to (127);
\draw[->] (127) to (128);

\draw[->] (134) to (145);
\draw[->] (145) to (156);
\draw[->] (156) to (167);
\draw[->] (167) to (178);

\draw[->] (145) to (124);
\draw[->] (156) to (125);
\draw[->] (167) to (126);

\draw[->] (124) to (134);
\draw[->] (125) to (145);
\draw[->] (126) to (156);
\draw[->] (127) to (167);

\draw[->] (134) to (234);
\draw[->] (145) to (345);
\draw[->] (156) to (456);
\draw[->] (167) to (567);

\draw[->] (178) to (127);
\draw[->] (678) to (167);
\draw[->] (567) to (156);
\draw[->] (456) to (145);
\draw[->] (345) to (134);

\node (GOctagon167) at (10.5,0) [regular polygon, regular polygon sides=8, draw,minimum size=0.6cm] {};
\fill (GOctagon167.corner 2) -- (GOctagon167.corner 5) -- (GOctagon167.corner 4);

\node (GOctagon127) [below=of GOctagon167,regular polygon, regular polygon sides=8, draw,minimum size=0.6cm] {};
\fill (GOctagon127.corner 2) -- (GOctagon127.corner 1) -- (GOctagon127.corner 4);

\node (GOctagon156) [right=of GOctagon167,regular polygon, regular polygon sides=8, draw,minimum size=0.6cm] {};
\fill (GOctagon156.corner 2) -- (GOctagon156.corner 6) -- (GOctagon156.corner 5);

\node (GOctagon126) [below=of GOctagon156,regular polygon, regular polygon sides=8, draw,minimum size=0.6cm] {};
\fill (GOctagon126.corner 2) -- (GOctagon126.corner 1) -- (GOctagon126.corner 5);

\node (GOctagon145) [right=of GOctagon156,regular polygon, regular polygon sides=8, draw,minimum size=0.6cm] {};
\fill (GOctagon145.corner 2) -- (GOctagon145.corner 7) -- (GOctagon145.corner 6);

\node (GOctagon125) [below=of GOctagon145,regular polygon, regular polygon sides=8, draw,minimum size=0.6cm] {};
\fill (GOctagon125.corner 2) -- (GOctagon125.corner 1) -- (GOctagon125.corner 6);

\node (GOctagon134) [right=of GOctagon145,regular polygon, regular polygon sides=8, draw,minimum size=0.6cm] {};
\fill (GOctagon134.corner 2) -- (GOctagon134.corner 8) -- (GOctagon134.corner 7);

\node (GOctagon124) [below=of GOctagon134,regular polygon, regular polygon sides=8, draw,minimum size=0.6cm] {};
\fill (GOctagon124.corner 2) -- (GOctagon124.corner 1) -- (GOctagon124.corner 7);

\node (GOctagon123) [right=of GOctagon124,regular polygon, regular polygon sides=8, draw,minimum size=0.6cm]  {};
\fill (GOctagon123.corner 2) -- (GOctagon123.corner 1) -- (GOctagon123.corner 8);

\node (GOctagon128) [left=of GOctagon127,regular polygon, regular polygon sides=8, draw,minimum size=0.6cm] {};
\fill (GOctagon128.corner 2) -- (GOctagon128.corner 1) -- (GOctagon128.corner 3);

\node (GOctagon178) [left=of GOctagon167,regular polygon, regular polygon sides=8, draw,minimum size=0.6cm] {};
\fill (GOctagon178.corner 2) -- (GOctagon178.corner 4) -- (GOctagon178.corner 3);

\node (GOctagon678) [above=of GOctagon178,regular polygon, regular polygon sides=8, draw,minimum size=0.6cm] {};
\fill (GOctagon678.corner 5) -- (GOctagon678.corner 4) -- (GOctagon678.corner 3);

\node (GOctagon567) [above=of GOctagon167,regular polygon, regular polygon sides=8, draw,minimum size=0.6cm] {};
\fill (GOctagon567.corner 6) -- (GOctagon567.corner 5) -- (GOctagon567.corner 4);

\node (GOctagon456) [above=of GOctagon156,regular polygon, regular polygon sides=8, draw,minimum size=0.6cm] {};
\fill (GOctagon456.corner 7) -- (GOctagon456.corner 6) -- (GOctagon456.corner 5);

\node (GOctagon345) [above=of GOctagon145,regular polygon, regular polygon sides=8, draw,minimum size=0.6cm] {};
\fill (GOctagon345.corner 8) -- (GOctagon345.corner 7) -- (GOctagon345.corner 6);

\node (GOctagon234) [above=of GOctagon134,regular polygon, regular polygon sides=8, draw,minimum size=0.6cm] {};
\fill (GOctagon234.corner 1) -- (GOctagon234.corner 8) -- (GOctagon234.corner 7);

\draw[->,shorten <= 1mm, shorten >= 1mm] (GOctagon123) to (GOctagon124);
\draw[->,shorten <= 1mm, shorten >= 1mm] (GOctagon124) to (GOctagon125);
\draw[->,shorten <= 1mm, shorten >= 1mm] (GOctagon125) to (GOctagon126);
\draw[->,shorten <= 1mm, shorten >= 1mm] (GOctagon126) to (GOctagon127);
\draw[->,shorten <= 1mm, shorten >= 1mm] (GOctagon127) to (GOctagon128);

\draw[->,shorten <= 1mm, shorten >= 1mm] (GOctagon134) to (GOctagon145);
\draw[->,shorten <= 1mm, shorten >= 1mm] (GOctagon145) to (GOctagon156);
\draw[->,shorten <= 1mm, shorten >= 1mm] (GOctagon156) to (GOctagon167);
\draw[->,shorten <= 1mm, shorten >= 1mm] (GOctagon167) to (GOctagon178);

\draw[->,shorten <= 1mm, shorten >= 1mm] (GOctagon145) to (GOctagon124);
\draw[->,shorten <= 1mm, shorten >= 1mm] (GOctagon156) to (GOctagon125);
\draw[->,shorten <= 1mm, shorten >= 1mm] (GOctagon167) to (GOctagon126);

\draw[->,shorten <= 1mm, shorten >= 1mm] (GOctagon124) to (GOctagon134);
\draw[->,shorten <= 1mm, shorten >= 1mm] (GOctagon125) to (GOctagon145);
\draw[->,shorten <= 1mm, shorten >= 1mm] (GOctagon126) to (GOctagon156);
\draw[->,shorten <= 1mm, shorten >= 1mm] (GOctagon127) to (GOctagon167);

\draw[->,shorten <= 1mm, shorten >= 1mm] (GOctagon134) to (GOctagon234);
\draw[->,shorten <= 1mm, shorten >= 1mm] (GOctagon145) to (GOctagon345);
\draw[->,shorten <= 1mm, shorten >= 1mm] (GOctagon156) to (GOctagon456);
\draw[->,shorten <= 1mm, shorten >= 1mm] (GOctagon167) to (GOctagon567);

\draw[->,shorten <= 1mm, shorten >= 1mm] (GOctagon178) to (GOctagon127);
\draw[->,shorten <= 1mm, shorten >= 1mm] (GOctagon678) to (GOctagon167);
\draw[->,shorten <= 1mm, shorten >= 1mm] (GOctagon567) to (GOctagon156);
\draw[->,shorten <= 1mm, shorten >= 1mm] (GOctagon456) to (GOctagon145);
\draw[->,shorten <= 1mm, shorten >= 1mm] (GOctagon345) to (GOctagon134);

\end{tikzpicture}}
\end{center}
This quiver is mutation-equivalent to the Dynkin diagram of type $E_{8}$.  We again omit the corresponding quasi-commutation matrix; the compatibility calculation gives the matrix $(2I_{8}\ 0_{8,8})$.  (Note that if we considered the analogous cluster for $n=9$, we would have more mutable cluster variables than coefficients.  This is related to the fact that $\complex[\mathrm{Gr}(3,9)]$ is of infinite type.)

We will again not give the full list of quantum cluster variables.  By computer-aided calculation, we have verified that the set of quantum cluster variables in this example consists of 
\begin{itemize}
\item all quantum Pl\"{u}cker coordinates;
\item the 56 elements $X_{q}^{abcdef}$ and $Y_{q}^{abcdef}$ defined as above, where $\{a,b,c,d,e,f\}$ is now a subset of $\{1,\ldots ,8 \}$ written in increasing order; and
\item 24 quantizations of cubic regular functions, $A_{q}(i)$ ($1\leq i \leq 8$) and $B_{q}(i,j)$ ($1\leq i \leq 8$, $0\leq j \leq 1$), listed in Table~\ref{Gr38cubics}.
\end{itemize}

\renewcommand{\arraystretch}{1.25}
\begin{table}
\begin{center}
\[ \begin{array}{r@{\,=\,}l@{\,}l@{\,-\,}l@{\,}l@{\,-\,}l@{\,}l}
A_{q}(1) & q^{-1} & \qminor{134}\qminor{258}\qminor{167} & q & \qminor{134}\qminor{678}\qminor{125} & q & \qminor{158}\qminor{234}\qminor{167} \\
A_{q}(2) & q^{-1} & \qminor{245}\qminor{136}\qminor{278} & q & \qminor{245}\qminor{178}\qminor{236} & q^{-1} & \qminor{126}\qminor{345}\qminor{278} \\
A_{q}(3) & q^{} & \qminor{356}\qminor{247}\qminor{138} & q^{} & \qminor{356}\qminor{128}\qminor{347} & q^{} & \qminor{237}\qminor{456}\qminor{138} \\
A_{q}(4) & q^{3} & \qminor{467}\qminor{358}\qminor{124} & q^{} & \qminor{467}\qminor{123}\qminor{458} & q^{3} & \qminor{348}\qminor{567}\qminor{124} \\
A_{q}(5) & q^{3} & \qminor{578}\qminor{146}\qminor{235} & q^{3} & \qminor{578}\qminor{234}\qminor{156} & q^{} & \qminor{145}\qminor{678}\qminor{235} \\
A_{q}(6) & q^{} & \qminor{168}\qminor{257}\qminor{346} & q^{} & \qminor{168}\qminor{345}\qminor{267} & q^{} & \qminor{256}\qminor{178}\qminor{346} \\
A_{q}(7) & q^{-1} & \qminor{127}\qminor{368}\qminor{457} & q^{-1} & \qminor{127}\qminor{456}\qminor{378} & q^{} & \qminor{367}\qminor{128}\qminor{457} \\
A_{q}(8) & q^{-1} & \qminor{238}\qminor{147}\qminor{568} & q^{} & \qminor{238}\qminor{567}\qminor{148} & q^{} & \qminor{478}\qminor{123}\qminor{568} \\[2em]
B_{q}(1,0) &  & \qminor{258}\qminor{134}\qminor{267} & q^{2} & \qminor{258}\qminor{167}\qminor{234} & q^{-3} & \qminor{128}\qminor{234}\qminor{567} \\
B_{q}(2,0) & q^{-2} & \qminor{136}\qminor{245}\qminor{378} &  & \qminor{136}\qminor{278}\qminor{345} & q^{-5} & \qminor{123}\qminor{345}\qminor{678} \\
B_{q}(3,0) &  & \qminor{247}\qminor{356}\qminor{148} & & \qminor{247}\qminor{138}\qminor{456} & q^{-3} & \qminor{234}\qminor{456}\qminor{178} \\
B_{q}(4,0) & q^{2} & \qminor{358}\qminor{467}\qminor{125} &  & \qminor{358}\qminor{124}\qminor{567} & q^{-1} & \qminor{345}\qminor{567}\qminor{128} \\
B_{q}(5,0) &  & \qminor{146}\qminor{578}\qminor{236} & q^{-2} & \qminor{146}\qminor{235}\qminor{678} & q^{} & \qminor{456}\qminor{678}\qminor{123} \\
B_{q}(6,0) &  & \qminor{257}\qminor{168}\qminor{347} & & \qminor{257}\qminor{346}\qminor{178} & q^{} & \qminor{567}\qminor{178}\qminor{234} \\
B_{q}(7,0) &  & \qminor{368}\qminor{127}\qminor{458} & q^{2} & \qminor{368}\qminor{457}\qminor{128} & q^{} & \qminor{678}\qminor{128}\qminor{345} \\
B_{q}(8,0) & & \qminor{147}\qminor{238}\qminor{156} & q^{2} & \qminor{147}\qminor{568}\qminor{123} & q^{-1} & \qminor{178}\qminor{123}\qminor{456} \\[2em]
B_{q}(1,1) & & \qminor{258}\qminor{167}\qminor{348} & q^{-2} & \qminor{258}\qminor{134}\qminor{678} & q^{} & \qminor{128}\qminor{678}\qminor{345} \\
B_{q}(2,1) &  & \qminor{136}\qminor{278}\qminor{145} & q^{-2} & \qminor{136}\qminor{245}\qminor{178} & q^{-1} & \qminor{123}\qminor{178}\qminor{456} \\
B_{q}(3,1) &  & \qminor{247}\qminor{138}\qminor{256} &  & \qminor{247}\qminor{356}\qminor{128} & q^{-1} & \qminor{234}\qminor{128}\qminor{567} \\
B_{q}(4,1) &  & \qminor{358}\qminor{124}\qminor{367} & q^{2} & \qminor{358}\qminor{467}\qminor{123} & q^{-1} & \qminor{345}\qminor{123}\qminor{678} \\
B_{q}(5,1) & q^{-2} & \qminor{146}\qminor{235}\qminor{478} & & \qminor{146}\qminor{578}\qminor{234} & q^{} & \qminor{456}\qminor{234}\qminor{178} \\
B_{q}(6,1) &  & \qminor{257}\qminor{346}\qminor{158} & & \qminor{257}\qminor{168}\qminor{345} & q^{3} & \qminor{567}\qminor{345}\qminor{128} \\
B_{q}(7,1) & q^{2} & \qminor{368}\qminor{457}\qminor{126} & & \qminor{368}\qminor{127}\qminor{456} & q^{5} & \qminor{678}\qminor{456}\qminor{123} \\
B_{q}(8,1) & & \qminor{147}\qminor{568}\qminor{237} & q^{-2} & \qminor{147}\qminor{238}\qminor{567} & q^{3} & \qminor{178}\qminor{567}\qminor{234}
\end{array} \]
\caption{The quantizations of cubic regular functions occurring as quantum cluster variables for $\mathbb{K}_{q}[\mathrm{Gr}(3,8)]$}\label{Gr38cubics}
\end{center}
\end{table}
\renewcommand{\arraystretch}{1}

We conclude that $\mathbb{K}_{q}[\mathrm{Gr}(3,8)]$ is a quantum cluster algebra of type $E_{8}$, provided that $q$ has a square root in $\mathbb{K}$.

We have observed a quantum analogue of the action of the dihedral group $D_{2n}$ on the affine cone $X(k,n)$ of the Grassmannian.  In the case of $\complex[\mathrm{Gr}(3,8)]$, Scott (\cite{Scott-Grassmannians}) notes that from knowledge of two cubic regular functions $A$ and $B$, one may obtain the remaining twenty-two cubic functions that are also cluster variables by applying to their indices the permutations $\rho=(1\, 8\, 7\, 6\, 5\, 4\, 3\, 2)$ and $\sigma=(2\, 8)(3\, 7)(4\, 6)$.  This gives two families: eight cluster variables $A^{\rho^{r}}$, $0\leq r \leq 7$ and sixteen cluster variables $B^{\sigma^{r}\rho^{s}}$, $0\leq r \leq 1$, $0\leq s \leq 7$ ($A$ is invariant under $\sigma$).

As noted in \cite{LL-twistingGr}, index-cycling is not an automorphism of the quantum Grassmannian but there is a cocycle twist that replaces it.  From the explicit calculations we have observed that one passes from the quantum cluster variable $A_{q}(i)$ to $A_{q}(i+1)$ (respectively, $B_{q}(i,0)$ to $B_{q}(i+1,0)$ and $B_{q}(i,1)$ to $B_{q}(i+1,1)$) precisely by means of this cocycle twist up to powers of $q$.  (The twist raises indices, thus corresponds to $\rho^{-1}$.)

We would then also expect a quantum version of $\sigma$ to take us between $B_{q}(i,0)$ and $B_{q}(i,1)$.  It is perhaps not surprising that we have been able to find an expression for $B_{q}(1,1)$ with quantum minor indices being those of $B_{q}(1,0)$ under the action of $\sigma$, but with only this one expression to work with we cannot yet conclude that there is a suitable cocycle to quantize $\sigma$.  Since the dihedral action exists classically for all Grassmannians, it ought to be possible to see its quantum analogue in other cases than just $\mathbb{K}_{q}[\mathrm{Gr}(3,8)]$ and so this remains a topic for further investigation.

\section{Quantum Schubert cells}\label{s:qSchubertCells}

In \cite{LLR}, motivated by the classical setting, the notion of quantum Schubert cells in the quantum Grassmannian was defined via noncommutative dehomogenisation as follows. First recall that the standard order on the set of quantum Pl\" ucker coordinates of $\mathbb{K}_q[\mathrm{Gr}(2,n)]$ is defined by
\[ \qminor{ij} \leq \qminor{kl} \mbox{ when } i \leq k \mbox{ and } j \leq l.\] 
Fix a  quantum Pl\" ucker coordinate $\delta$ of $\mathbb{K}_q[\mathrm{Gr}(2,n)]$. The quantum Schubert ideal $I_{\delta}$ corresponding to $\delta$ is the two-sided ideal of $\mathbb{K}_q[\mathrm{Gr}(2,n)]$ generated by the quantum Pl\" ucker coordinates $\gamma$ with $\gamma \ngeq \delta$. In the quantum Schubert variety $\mathbb{K}_q[\mathrm{Gr}(2,n)]/ I_{\delta}$, the image $\overline{\delta}$ of $\delta$ in this factor algebra is normal and regular, so that one can form the localisation 
 \[S_{\delta} \defeq \frac{\mathbb{K}_q[\mathrm{Gr}(2,n)]}{I_{\delta}}\left[ \overline{\delta}^{-1} \right]. \]
Now observe that $\mathbb{K}_q[\mathrm{Gr}(2,n)]$ is graded in a natural way with all quantum Pl\" ucker coordinates in degree 1. As $I_{\delta}$ and $\delta$ are homogeneous, $S_{\delta}$ is a $\mathbb{Z}$-graded algebra. Then the quantum Schubert cell $S^0(\delta)$ associated to $\delta$ is the degree $0$ part of $S_{\delta}$ \cite[Definition 4.1]{LLR}. 

It was proved in \cite{LLR} that this algebra can be identified to a subalgebra of $\mathbb{K}_{q^{-1}}[\mathrm{M}(2,n-2)]$.  Set $\delta=\qminor{kl}$. Then \cite[Theorem 4.7]{LLR} shows that $S^0(\delta)$ is isomorphic to the subalgebra of $ \mathbb{K}_{q^{-1}}[\mathrm{M}(2,n-2)]$ generated by $X_{1j}$ with $1 \leq j \leq n-1-k$ and $X_{2j}$ with $1 \leq j \leq n-l$. We refer the reader to \cite{LLR} for more details.

Of course, we can interchange $q$ and $q^{-1}$, so that it is natural to refer to the subalgebra of 
$ \mathbb{K}_{q}[\mathrm{M}(2,n-2)]$ generated by $X_{1j}$ with $1 \leq j \leq t$ and $X_{2j}$ with $1 \leq j \leq s$ as the 
quantum Schubert cell of the Grassmannian $\mathrm{Gr}(2,n)$ associated to the partition $(t,s)$, where $t\geq s$ and  $s,t\leq n-2$. 
Our aim is to prove that this algebra is a quantum cluster algebra.

 First, we state our chosen initial seed for the quantum cluster algebra structure.  For our initial quantum cluster we choose
\[ \underline{y} = ( X_{1t},X_{1(t-1)},\ldots,X_{1s},\matentrybf{1(\itbf{s}-1)},\matentrybf{1(\itbf{s}-2)},\ldots,\matentrybf{11},X_{21},\qminor{12},\qminor{23},\ldots,\qminor{(s-1)s}). \]
So, we have $s-1$ mutable cluster variables and $t+1$ coefficients (and a cluster containing $s+t$ variables in total).

We will retain a similar notation to the Grassmannian case, except that rows and columns corresponding to $X_{ij}$ will be indicated by parenthesis $(ij)$ and those to minors by brackets $[ij]$, as before.  The initial exchange matrix $B$ is given by
\begin{align*}
B_{(1i)(1k)} & =  
  \begin{cases} 
    0 & \text{for}\ i\geq k+2 \\
   -1 & \text{for}\ i=k+1 \\
    0 & \text{for}\ i=k \\
    1 & \text{for}\ i=k-1 \\
    0 & \text{for}\ i\leq k-2
  \end{cases} &
B_{[j(j+1)](1k)} & =
  \begin{cases} 
    0 & \text{for}\ j\leq k-2 \\
   -1 & \text{for}\ j=k-1 \\
    1 & \text{for}\ j=k \\
    0 & \text{for}\ j\geq k+1
  \end{cases} \\
B_{(21)(1k)} & =
  \begin{cases}
    0 & \text{for}\ k\geq 2 \\
   -1 & \text{for}\ k=1
  \end{cases} & &
\end{align*}
where $t\geq i\geq 1$, $1\leq j \leq s-1$ and $s-1\geq k\geq 1$.  For example, for $n\geq 9$, $t=7$ and $s=6$ we have 
\[ B = 
\bordermatrix{
 & \scriptstyle{(15)} & \scriptstyle{(14)} & \scriptstyle{(13)} & \scriptstyle{(12)} & \scriptstyle{(11)} \cr
\scriptstyle{(17)} & 0 & 0 & 0 & 0 & 0 \cr
\scriptstyle{(16)} & -1 & 0 & 0 & 0 & 0 \cr
\scriptstyle{(15)} & \mathbf{0} & \mathbf{-1} & \mathbf{0} & \mathbf{0} & \mathbf{0} \cr
\scriptstyle{(14)} & \mathbf{1} & \mathbf{0} & \mathbf{-1} & \mathbf{0} & \mathbf{0}  \cr
\scriptstyle{(13)} & \mathbf{0} & \mathbf{1} & \mathbf{0} & \mathbf{-1} & \mathbf{0}  \cr
\scriptstyle{(12)} & \mathbf{0} & \mathbf{0} & \mathbf{1} & \mathbf{0} & \mathbf{-1}  \cr
\scriptstyle{(11)} & \mathbf{0} & \mathbf{0} & \mathbf{0} & \mathbf{1} & \mathbf{0}  \cr
\scriptstyle{(21)} & 0 & 0 & 0 & 0 & -1 \cr
\scriptstyle{[12]} & 0 & 0 & 0 & -1 & 1 \cr
\scriptstyle{[23]} & 0 & 0 & -1 & 1 & 0 \cr
\scriptstyle{[34]} & 0 & -1 & 1 & 0 & 0 \cr
\scriptstyle{[45]} & -1 & 1 & 0 & 0 & 0 \cr
\scriptstyle{[56]} & 1 & 0 & 0 & 0 & 0 
} \]
We note that $B$ has a natural three-block structure, as a zero block $0_{t-s,s-1}$ on the row set $\{(1t),\ldots,(1(s+1))\}$ and two blocks on the row sets $\{(1s),\ldots,(11)\}$ and $\{(21),[12],\ldots,[(s-1)s]\}$.

The corresponding quiver for $n\geq 9$, $t=7$ and $s=6$ is 
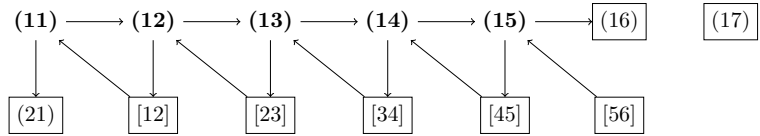
\begin{center} 
\scalebox{0.75}{\begin{tikzpicture}

\node (x11) at (0,2) {\textbf{(11)}};
\node (x12) [right=of x11] {\textbf{(12)}};
\node (x13) [right=of x12] {\textbf{(13)}};
\node (x14) [right=of x13] {\textbf{(14)}};
\node (x15) [right=of x14] {\textbf{(15)}};
\node (x16) [rectangle,draw=black] [right=of x15] {(16)};
\node (x17) [rectangle,draw=black] [right=of x16] {(17)};
\node (x21) [rectangle,draw=black] [below=of x11] {(21)};
\node (m12) [rectangle,draw=black] [below=of x12] {[12]};
\node (m23) [rectangle,draw=black] [below=of x13] {[23]};
\node (m34) [rectangle,draw=black] [below=of x14] {[34]};
\node (m45) [rectangle,draw=black] [below=of x15] {[45]};
\node (m56) [rectangle,draw=black] [below=of x16] {[56]};

\draw[->] (x11) to (x12);
\draw[->] (x11) to (x21);
\draw[->] (x12) to (x13);
\draw[->] (x12) to (m12);
\draw[->] (x13) to (x14);
\draw[->] (x13) to (m23);
\draw[->] (x14) to (x15);
\draw[->] (x14) to (m34);
\draw[->] (x15) to (x16);
\draw[->] (x15) to (m45);

\draw[->] (m12) to (x11);
\draw[->] (m23) to (x12);
\draw[->] (m34) to (x13);
\draw[->] (m45) to (x14);
\draw[->] (m56) to (x15);

\end{tikzpicture}}
\end{center}
We see that for any $n$ and partition $(t,s)$ this quantum cluster algebra is of type $A_{s-1}$ (independent of $t$), since the subquiver on the vertices $\{\mathbf{(11)},\ldots,\mathbf{(1(\itbf{s}-1))}\}$ is an orientation of the Dynkin diagram of this type.

The quasi-commutation matrix $L$ is as follows:
\allowdisplaybreaks
\begin{align*}
  L_{(1i)(1k)} & =
     \begin{cases}
       -1 & \text{for}\ i\geq k+1 \\
        0 & \text{for}\ i=k \\
        1 & \text{for}\ i\leq k-1
     \end{cases} &
  L_{(1i)(21)} & = 
     \begin{cases}
        0 & \text{for}\ i\geq 2 \\ 
        1 & \text{for}\ i=1
     \end{cases} \\
  L_{(1i)[l(l+1)]} & =
     \begin{cases}
       -1 & \text{for}\ i\geq l+2 \\
        0 & \text{for}\ i=l\ \text{or}\ l+1 \\
        1 & \text{for}\ i\leq l-1
     \end{cases} &
  L_{(21)(1k)} & =
     \begin{cases}
        0 & \text{for}\ k\geq 2 \\
       -1 & \text{for}\ k=1
     \end{cases} \\
  L_{(21)(21)} & = 0 &
  L_{(21)[l(l+1)]} & = 
     \begin{cases}
        0 & \text{for}\ l=1 \\
        1 & \text{for}\ l\geq 2
     \end{cases} \\
  L_{[j(j+1)](1k)} & =
     \begin{cases}
        1 & \text{for}\ j\leq k-2 \\
        0 & \text{for}\ j=k-1\ \text{or}\ k \\
       -1 & \text{for}\ j\geq k+1
     \end{cases} &
  L_{[j(j+1)](21)} & =
     \begin{cases}
        0 & \text{for}\ j=1 \\
       -1 & \text{for}\ j\geq 2
     \end{cases} \\
  L_{[j(j+1)][l(l+1)]} & =
     \begin{cases}
        2 & \text{for}\ j\leq l-2 \\
        1 & \text{for}\ j=l-1 \\
        0 & \text{for}\ j=l \\
       -1 & \text{for}\ j=l+1 \\
       -2 & \text{for}\ j\geq l+2
     \end{cases} &
\end{align*}
for $t\geq i,k\geq 1$ and $1\leq j,l \leq s-1$.  These values may be verified easily, using the well-known quasi-commutation relations in quantum matrices.  For $n \geq 9$ and the partition $(7,6)$, this matrix is
\[ L = 
\bordermatrix{
 & \scriptstyle{(17)} & \scriptstyle{(16)} & \scriptstyle{(15)} & \scriptstyle{(14)} & \scriptstyle{(13)} & \scriptstyle{(12)} & \scriptstyle{(11)} & \scriptstyle{(21)} & \scriptstyle{[12]} & \scriptstyle{[23]} & \scriptstyle{[34]} & \scriptstyle{[45]} & \scriptstyle{[56]} \cr
\scriptstyle{(17)} & 0 & -1 & -1 & -1 & -1 & -1 & -1 & 0 & -1 & -1 & -1 & -1 & -1 \cr
\scriptstyle{(16)} & 1 & 0 & -1 & -1 & -1 & -1 & -1 & 0 & -1 & -1 & -1 & -1 & 0 \cr
\scriptstyle{(15)} & 1 & 1 & \mathbf{0} & \mathbf{-1} & \mathbf{-1} & \mathbf{-1} & \mathbf{-1} & 0 & -1 & -1 & -1 & 0 & 0 \cr
\scriptstyle{(14)} & 1 & 1 & \mathbf{1} & \mathbf{0} & \mathbf{-1} & \mathbf{-1} & \mathbf{-1} & 0 & -1 & -1 & 0 & 0 & 1 \cr
\scriptstyle{(13)} & 1 & 1 & \mathbf{1} & \mathbf{1} & \mathbf{0} & \mathbf{-1} & \mathbf{-1} & 0 & -1 & 0 & 0 & 1 & 1 \cr
\scriptstyle{(12)} & 1 & 1 & \mathbf{1} & \mathbf{1} & \mathbf{1} & \mathbf{0} & \mathbf{-1} & 0 & 0 & 0 & 1 & 1 & 1 \cr
\scriptstyle{(11)} & 1 & 1 & \mathbf{1} & \mathbf{1} & \mathbf{1} & \mathbf{1} & \mathbf{0} &  1 &  0 &  1 &  1 &  1 & 1 \cr
\scriptstyle{(21)} &  0 &  0 &  0 &  0 &  0 &  0 & -1 &  0 &  0 &  1 &  1 &  1 & 1 \cr
\scriptstyle{[12]} &  1 &  1 &  1 &  1 &  1 &  0 & 0 &  0 &  0 &  1 &  2 &  2 & 2 \cr
\scriptstyle{[23]} &  1 &  1 &  1 &  1 &  0 &  0 & -1 & -1 & -1 &  0 &  1 &  2 & 2 \cr
\scriptstyle{[34]} &  1 &  1 &  1 &  0 &  0 & -1 & -1 & -1 & -2 & -1 &  0 &  1 & 2 \cr
\scriptstyle{[45]} &  1 &  1 &  0 &  0 & -1 & -1 & -1 & -1 & -2 & -2 & -1 &  0 & 1 \cr
\scriptstyle{[56]} &  1 &  0 &  0 & -1 & -1 & -1 & -1 & -1 & -2 & -2 & -2 & -1 & 0 
} \]

We claim that $B$ and $L$ are compatible.

\begin{proposition} $B^{T}L=(0_{s-1,t-s+1}\ 2I_{s-1}\ 0_{s-1,s})$
\end{proposition}

\begin{proof}  This follows from a similar set of calculations to those in the proof of Proposition~\ref{Gr-mats-compat}.
\end{proof}

To show that the quantum Schubert cells are quantum cluster algebras, we use a similar strategy to the Grassmannian case.  Firstly, we note that we can reduce to considering the case $t=s$.  For if $t>s$, we have additional coefficients $X_{1t},\ldots,X_{1(s+1)}$ but the corresponding rows of the exchange matrix $B$ are all zero and so these extra coefficients never appear in exchange relations.  Secondly, we can further reduce to the case $t=s=n-2$ (\ie the big cell), else we may as well consider a smaller $n$.

So, taking $t=s=n-2$, we are considering the initial quantum cluster
\[ \underline{y} = (X_{1(n-2)},\matentrybf{1(\itbf{n}-3)},\matentrybf{1(\itbf{n}-4)},\ldots,\matentrybf{11},X_{21},\qminor{12},\qminor{23},\ldots,\qminor{(n-3)(n-2)}). \]
We remark that the exchange matrix $B$ for this cluster is equal to that for the initial seed of $\mathbb{K}_{q}[\mathrm{Gr}(2,n)]$ with the row labelled $[12]$ deleted.  This is not surprising, as $\mathbb{K}_{q}[\mathrm{Gr}(2,n)][(\qminor{12})^{-1}]$ is isomorphic to a skew-Laurent extension of the form $\mathbb{K}_{q}[\mathrm{M}(2,n-2)][y^{\pm};\sigma]$ by \cite[Corollary 4.4]{KLR}. (The reader is referred to \cite{KLR} for more details on this isomorphism.)

We break the proof into two parts.  In the first, we show that we can obtain all remaining matrix entries by repeated mutation starting from $\matentrybf{11}$ and in the second we obtain the remaining minors.

\begin{proposition} Let $1\leq r\leq n-3$.  After the sequence of mutations $\mu_{(1r)}\circ \mu_{(1(r-1))}\circ \cdots \circ \mu_{(11)}$, the cluster obtained contains the mutable variables 
\[ (\matentrybf{1(\itbf{n}-3)},\ldots,\matentrybf{1(\itbf{r}+1)},\matentrybf{2(\itbf{r}+1)},\matentrybf{2\itbf{r}},\ldots,\matentrybf{22}).\]
Furthermore, the exchange matrix $B$ of this cluster has only the following non-zero entries in its $(1(r+1))$-column: $B_{(1(r+2))(1(r+1))}=B_{(2(r+1))(1(r+1))}=-1$ and $B_{[(r+1)(r+2)](1(r+1))}=1$.
\end{proposition}

\begin{proof}
We work by induction on $r$ and compute the mutation $\mu_{(1(r+1))}$ applied to the cluster and exchange matrix in the statement above.  This mutation exchanges $\matentrybf{(1(\itbf{r}+1))}$ for some new element $X$ which is determined by the quantum exchange relation.  This relation is determined by the non-zero entries of $B$ stated above and so is
\begin{eqnarray*} 
X & = & M\!\!\!\bordermatrix{ & & & \scriptscriptstyle{(1(r+1))} & & &  & \scriptscriptstyle{[(r+1)(r+2)]} & & \ \cr & \cdots & 0 & -1 & 0 & \cdots & 0 & 1 & 0 & \cdots} \\ 
  &  & +M\!\!\!\bordermatrix{ & & & \scriptscriptstyle{(1(r+2))} & \scriptscriptstyle{(1(r+1))} & & & & \scriptscriptstyle{(2(r+1))} & & \cr & \cdots & 0 & 1 & -1 & 0 & \cdots & 0 & 1 & 0 & \cdots } \\
  & = & q^{\alpha}\matentrybf{1(\itbf{r}+1)}^{-1}\qminorbf{(\itbf{r}+1)(\itbf{r}+2)}+q^{\beta}\matentrybf{1(\itbf{r}+2)}\matentrybf{1(\itbf{r}+1)}^{-1}\matentrybf{2(\itbf{r}+1)}
\end{eqnarray*}
where $\alpha$ and $\beta$ are integers to be calculated from $B$ and the quasi-commutation matrix.  These integers have the following values:
\begin{align*} \alpha &  =  \textstyle{\frac{1}{2}}((-1)\cdot (1)\cdot L_{[(r+1)(r+2)](1(r+1))})  = \textstyle{\frac{1}{2}}(0) =  0 \\
\intertext{and}
\beta & = \textstyle{\frac{1}{2}}(1\cdot (-1)\cdot L_{(1(r+1))(1(r+2))}+1\cdot 1\cdot L_{(2(r+1)(1(r+2))}+(-1)\cdot 1\cdot L_{(2(r+1))(1(r+1))}) \\
 & = \textstyle{\frac{1}{2}}(-1+0+1) \\
 & = 0
\end{align*}
(remembering that the $X_{ij}$ satisfy the quantum matrix relations).
Substituting these and applying quasi-commutation relations to move the $\matentrybf{1(\itbf{r}+1)}^{-1}$ terms to the left before re-arranging, we obtain the following equality involving $X$:
\[ \matentrybf{1(\itbf{r}+1)}X = \qminorbf{(\itbf{r}+1)(\itbf{r}+2)}+q\matentrybf{1(\itbf{r}+2)}\matentrybf{2(\itbf{r}+1)}. \]
But the right-hand side of this equation is equal to $\matentrybf{1(\itbf{r}+1)}\matentrybf{2(\itbf{r}+2)}$ and hence we deduce that $X=\matentrybf{2(\itbf{r}+2)}$ as expected.

It remains to verify the form of the $(1(r+2))$ column of the exchange matrix $B'$ after doing $\mu_{(1(r+1))}$.  Prior to this mutation, the $(1(r+2))$ column has zero entries except in rows $(1(r+3))$ and $[(r+1)(r+2)]$ where it is $-1$ and in rows $(2(r+2))$ and $[(r+2)(r+3)]$ where it is $1$.  (The $(2(r+2))$ label replaces $(1(r+1))$.)  The $(2(r+1))$-column is the $(1(r+1))$-column with signs changed and it is the unchanged one that we need to do the mutation: this has a $-1$ in rows $(1(r+2))$ and $(2(r+1))$ and a $1$ in row $[(r+1)(r+2)]$.  Therefore, after mutation, the $(1(r+2))$-column is zero except possibly in the following entries:
\begin{enumerate}[label=(\roman*)]
\item $-1$ in row $(1(r+3))$, unchanged since there is a $0$ in the $(1(r+1))$-column 
\item $0$ in row $(1(r+2))$, unchanged since the $(1(r+1))$-row and column values are $1$ and $-1$ respectively
\item $-1$ in row $(2(r+2))$, changed from $1$ by the sign change in the mutation row
\item $0$ in row $(2(r+1))$,  unchanged since the $(1(r+1))$-row and column values are $1$ and $-1$ respectively
\item $0$ in row $[(r+1)(r+2)]$, changed from $-1$ since the $(1(r+1))$-row and column values are both $1$
\item $1$ in row $[(r+2)(r+3)]$, unchanged since there is a $0$ in the $(1(r+1))$-column 
\end{enumerate}
In summary, the non-zero entries in column $(1(r+2))$ are $B'_{(1(r+3))(1(r+2))}=B'_{(2(r+2))(1(r+2))}=-1$ and $B'_{[(r+2)(r+3)](1(r+2))}=1$, verifying the inductive claim.
\end{proof}

Next, we state and prove the corresponding result on the minors.

\begin{proposition} Let $2\leq \alpha<\beta\leq n-3$.  After the sequence of mutations $\mu_{(1\beta)}\circ \mu_{(1(\beta-1))}\circ \cdots \circ \mu_{(1\alpha)}$, the cluster obtained contains the mutable variables 
\[ (\matentrybf{1(\itbf{n}-3)},\ldots,\matentrybf{1(\bbeta+1)},\qminorbf{(\balpha-1)(\bbeta+1)},\qminorbf{(\balpha-1)\bbeta},\ldots,\qminorbf{(\balpha-1)(\balpha+1)},\matentrybf{1(\balpha-1)},\ldots,\matentrybf{11}). \]
Furthermore, the exchange matrix $B$ of this cluster has only the following non-zero entries in its $(1(\beta+1))$-column: $B_{(1(\beta+2))(1(\beta+1))}=B_{[(\alpha-1)(\beta+1)](1(\beta+1))}=-1$ and $B_{(1(\alpha-1))(1(\beta+1))}=B_{[(\beta+1)(\beta+2)](1(\beta+1))}=1$.
\end{proposition}

\begin{proof}
We work by induction on $\beta$ and compute the mutation $\mu_{(1(\beta+1))}$ applied to the cluster and exchange matrix in the statement above.  This mutation exchanges $\matentrybf{(1(\bbeta+1))}$ for some new element $X$ which is determined by the quantum exchange relation.  This relation is determined by the non-zero entries of $B$ stated above and so is
\begin{eqnarray*} 
X & = & M\!\!\!\bordermatrix{ & & & \scriptscriptstyle{(1(\beta+1))} & & &  & \scriptscriptstyle{(1(\alpha-1))} & & & & \scriptscriptstyle{[(\beta+1)(\beta+2)]} & & \ \cr & \cdots & 0 & -1 & 0 & \cdots & 0 & 1 & 0 & \cdots & 0 & 1 & 0 & \cdots } \\ 
  &  & +M\!\!\!\bordermatrix{ & & & \scriptscriptstyle{(1(\beta+2))} & \scriptscriptstyle{(1(\beta+1))} & \scriptscriptstyle{[(\alpha-1)(\beta+1)]} & & \cr & \cdots & 0 & 1 & -1 & 1 & 0 & \cdots } \\
  & = & q^{r}\matentrybf{1(\bbeta+1)}^{-1}\matentrybf{1(\balpha-1)}\qminorbf{(\bbeta+1)(\bbeta+2)}+q^{s}\matentrybf{1(\bbeta+2)}\matentrybf{1(\bbeta+1)}^{-1}\qminorbf{(\balpha-1)(\bbeta+1)}
\end{eqnarray*}
where $r$ and $s$ are integers to be calculated from $B$ and the quasi-commutation matrix.  These integers have the following values:
\begin{align*} r & = \textstyle{\frac{1}{2}}((-1)\cdot 1\cdot L_{(1(\alpha-1))(1(\beta+1))}+(-1)\cdot 1\cdot L_{[(\beta+1)(\beta+2)](1(\beta+1))}+1\cdot 1\cdot L_{[(\beta+1)(\beta+2)](1(\alpha-1))}) \\ & = \textstyle{\frac{1}{2}}(-1+0-1) \\ & = -1 \\
\intertext{and}
s & = \textstyle{\frac{1}{2}}(1\cdot (-1)\cdot L_{(1(\beta+1))(1(\beta+2))}+1\cdot 1\cdot L_{[(\alpha-1)(\beta+1)](1(\beta+2))}+(-1)\cdot 1\cdot L_{[(\alpha-1)(\beta+1)](1(\beta+1))}) \\
 & = \textstyle{\frac{1}{2}}(-1+1+0) \\
 & = 0
\end{align*}
(remembering that the $X_{ij}$ satisfy the quantum matrix relations and quasi-commutation relations involving minors can be calculated from these).

Substituting these and applying quasi-commutation relations to move the $\matentrybf{1(\bbeta+1)}^{-1}$ terms to the left before re-arranging, we obtain the following equality involving $X$:
\[ \matentrybf{1(\bbeta+1)}X = q^{-1}\matentrybf{1(\balpha-1)}\qminorbf{(\bbeta+1)(\bbeta+2)}+q\matentrybf{1(\bbeta+2)}\qminorbf{(\balpha-1)(\bbeta+1)}. \]
But the right-hand side of this equation is equal to $\matentrybf{1(\bbeta+1)}\qminorbf{(\balpha-1)(\bbeta+2)}$ by the $q$-Laplace relations (see for instance \cite[Lemma A.4]{GooLen}) and hence we deduce that $X=\qminorbf{(\balpha-1)(\bbeta+2)}$ as expected.

It remains to verify the form of the $(1(\beta+2))$ column of the exchange matrix $B'$ after doing $\mu_{(1(\beta+1))}$.  Prior to this mutation, the $(1(\beta+2))$ column has zero entries except in rows $(1(\beta+3))$ and $[(\alpha-1)(\beta+1)]$ where it is $-1$ and in rows $(1(\alpha-1))$ and $[(\beta+1)(\beta+2)]$ where it is $1$, by the inductive hypothesis.  (The $[(\alpha-1)(\beta+1)]$ label replaces $(1(\beta+1))$.)  The $[(\alpha-1)(\beta+1)]$-column is the $(1(\beta+1))$-column with signs changed and it is the unchanged one that we need to do the mutation: this has a $-1$ in rows $(1(\beta+2))$ and $[\beta(\beta+1)]$ and a $1$ in rows $(1\beta)$ and $[(\beta+1)(\beta+2)]$.  Therefore, after mutation, the $(1(\beta+2))$-column is zero except possibly in the following entries:
\begin{enumerate}[label=(\roman*)]
\item $-1$ in row $(1(\beta+3))$, unchanged since there is a $0$ in the $(1(\beta+1))$-column 
\item $0$ in row $(1(\beta+2))$, unchanged since the $(1(\beta+1))$-row and column values are $1$ and $-1$ respectively
\item $-1$ in row $[(\alpha-1)(\beta+2)]$, changed from $1$ by the sign change in the mutation row
\item $0$ in row $[(\alpha-1)(\beta+1)]$, unchanged since the $(1(\beta+1))$-row and column values are $1$ and $-1$ respectively
\item $1$ in row $(1(\alpha-1))$, changed from $0$ since the $(1(\beta+1))$-row and column values are both $1$
\item $0$ in row $[(\beta+1)(\beta+2)]$, changed from $1$ since the $(1(\beta+1))$-row and column values are both $1$
\item $1$ in row $[(\beta+2)(\beta+3)]$, unchanged since there is a $0$ in the $(1(\beta+1))$-column 
\end{enumerate}
That is, the non-zero entries in column $(1(\beta+2))$ are $B'_{(1(\beta+3))(1(\beta+2))}=B'_{[(\alpha-1)(\beta+2)](1(\beta+2))}=-1$ and $B'_{(1(\alpha-1))(1(\beta+2))}=B'_{[(\beta+2)(\beta+3)](1(\beta+2))}=1$, verifying the inductive claim.
\end{proof}

We can now complete the proof of our claim.

\begin{theorem} The quantum Schubert cell of the Grassmannian $\mathrm{Gr}(2,n)$ associated to the partition $(t,s)$, where $t\geq s$ and $t,s\leq n-2$, is a quantum cluster algebra of type $A_{s-1}$.
\end{theorem}

\begin{proof}  We recall from above that we were able to reduce to the case $t=s=n-2$.  It follows from the previous propositions that every generator occurs in some cluster and can be obtained by a finite sequence of mutations from our given initial seed.  Furthermore, the general theory of (quantum) cluster algebras of type $A$ tells us that a cluster algebra of type $A_{l}$ has precisely $\frac{l(l+1)}{2}+l=\frac{l^{2}+3l}{2}$ mutable cluster variables (as this is the number of almost positive roots in this type).  For $l=s-1=n-3$, this number is $\frac{(n-3)(n-2)}{2}+(n-3)=\frac{n^{2}-3n}{2}$.  

Now, the number of matrix entries and minors are $2(n-2)$ and $\frac{(n-3)(n-2)}{2}$ respectively and our stated cluster algebra structure has $n-1$ coefficients.  Since $2n-4+\frac{(n-3)(n-2)}{2}-(n-1)=\frac{n^{2}-3n}{2}$, we deduce that the set of all mutable cluster variables must be equal to the set of matrix entries and minors minus the coefficients, since the latter is contained in the former by the propositions.  

In other words, taking all cluster variables (mutable ones and coefficients) together gives us exactly the set of all matrix entries and minors.  From this we deduce that every cluster variable is a genuine matrix entry or minor (\ie no localisation is required) and hence the subalgebra generated by all cluster variables is exactly the whole quantum Schubert cell, which is therefore a quantum cluster algebra.
\end{proof}

\begin{remark} As for the quantum Grassmannian $\mathbb{K}_{q}[\mathrm{Gr}(2,n)]$ (see Remark \ref{qhypothesis}), the hypothesis that $q^{1/2}$ exists in $\mathbb{K}$ can be removed.
\end{remark}

\small

\bibliographystyle{amsplain}
\bibliography{biblio}

\normalsize

\end{document}